\documentclass[aos,seceqn,nameyear,dvips]{arximspdf}
\usepackage{mathrsfs}
\usepackage{graphicx}

\doi{10.1214/10-AOS813}
\volume{38}
\issue{6}
\pubyear{2010}
\firstpage{3321}
\lastpage{3351}

\makeatletter

\renewcommand{\mathcal}{\mathscr}

\newtheorem{theorem}{Theorem}[section]
\newtheorem{cor}[theorem]{Corollary}
\newtheorem{lm}{Lemma}

\newcommand{\iint}{\int\!\!\int}

\makeatother

\begin{document}
\begin{frontmatter}

\title{Uniform convergence rates for nonparametric regression and principal
component analysis in functional/longitudinal data}
\runtitle{Uniform convergence rates for functional data}

\begin{aug}
\author[A]{\fnms{Yehua} \snm{Li}\thanksref{t1}\ead[label=e1]{yehuali@uga.edu}} and
\author[B]{\fnms{Tailen} \snm{Hsing}\corref{}\thanksref{t2}\ead[label=e2]{thsing@umich.edu}}
\runauthor{Y. Li and T. Hsing}
\affiliation{University of Georgia and University of Michigan}
\address[A]{Department of Statistics\\
University of Georgia\\
Athens, Georgia 30602-7952 \\
USA\\
\printead{e1}}
\address[B]{Department of Statistics\\
University of Michigan\\
Ann Arbor, Michigan 48109-1107\\
USA\\
\printead{e2}}
\end{aug}

\thankstext{t1}{Supported by NSF Grant DMS-08-06131.}

\thankstext{t2}{Supported by NSF Grants DMS-08-08993 and DMS-08-06098.}

\received{\smonth{10} \syear{2009}}
\revised{\smonth{2} \syear{2010}}

%
\begin{abstract}
We consider nonparametric estimation of the mean and covariance
functions for functional/longitudinal data. Strong uniform convergence
rates are developed for estimators that are local-linear smooth\-ers. Our
results are obtained in a unified framework in which the number of
observations within each curve/cluster can be of any rate relative to
the sample size. We show that the convergence rates for the procedures
depend on both the number of sample curves and the number of
observations on each curve. For sparse functional data, these rates are
equivalent to the optimal rates in nonparametric regression. For dense
functional data, root-$n$ rates of convergence can be achieved with
proper choices of bandwidths. We further derive almost sure rates of
convergence for principal component analysis using the estimated
covariance function. The results are illustrated with simulation
studies.
\end{abstract}

%
\begin{keyword}[class=AMS]
\kwd[Primary ]{62J05}
\kwd[; secondary ]{62G20}
\kwd{62M20}.
\end{keyword}
\begin{keyword}
\kwd{Almost sure convergence}
\kwd{functional data analysis}
\kwd{kernel}
\kwd{local polynomial}
\kwd{nonparametric inference}
\kwd{principal components}.
\end{keyword}

\end{frontmatter}

\section{Introduction}
Estimating the mean and covariance functions are essential problems in
longitudinal and functional data analysis. Many recent papers focused
on nonparametric estimation so as to model the mean and covariance
structures flexibly. A partial list of such work includes Ramsay and
Silverman (\citeyear{RS05}), \citet{LC00},
\citet{W03}, Yao, M\"uller and Wang (\citeyear{YMW05a},
\citeyear{YMW05b}), Yao and Lee (\citeyear{YL06}) and
\citet{HMW06}.

On the other hand, functional principal component analysis (FPCA) based
on nonparametric covariance estimation has become one of the most
common dimension reduction approaches in functional data analysis.
Applications include temporal trajectory interpolation [Yao, M\"uller
and Wang (\citeyear{YMW05a})], functional generalized linear models
[\citet{MS05} and Yao, M\"uller and Wang (\citeyear{YMW05b})] and
functional sliced inverse regression [F\'erre and Yao
(\citeyear{FY05}), Li and Hsing (\citeyear{LH09})], to name a few. A
number of algorithms have been proposed for FPCA, some of which are
based on spline smoothing [\citet{JHS00}, \citet{ZHC08}] and
others based on kernel smoothing [Yao, M\"uller and Wang
(\citeyear{YMW05a}), \citet{HMW06}]. As usual,
large-sample theories can provide a basis for understanding the
properties of these estimators. So far, the asymptotic theories for
estimators based on kernel smoothing or local-polynomial smoothing are
better understood than those based on spline smoothing.

Some definitive theoretical findings on FPCA emerged in recent years.
In particular, \citet{HH06} proved various asymptotic expansions
for FPCA for densely recorded functional data, and \citet{HMW06}
established the optimal $L^2$ convergence rate for FPCA in the sparse
functional data setting. One of the most interesting findings in
\citet{HMW06} was that the estimated eigenfunctions, although
computed from an estimated two-dimensional surface, enjoy the
convergence rate of one-dimensional smoothers, and under favorable
conditions the estimated eigenvalues are root-$n$ consistent. In
contrast with the $L^2$ convergence rates of these nonparametric
estimators, less is known in term of uniform convergence rates. Yao,
M\"uller and Wang (\citeyear{YMW05a}) studied the uniform consistency
of the estimated mean, covariance and eigenfunctions, and demonstrated
that such uniform convergence properties are useful in many settings;
some other examples can also be found in \citet{Letal08}.

In classical nonparametric regression where observations are
independent, there are a number of well-known results concerning the
uniform convergence rates of kernel-based estimators. Those include
\citet{BR73}, \citet{HJS88} and H\"ardle (\citeyear{H89}).
More recently, \citet{CK03} extended some of those results to
local likelihood estimators and local estimating equations. However, as
remarked in Yao, M\"uller and Wang (\citeyear{YMW05a}), whether those
optimal rates can be extended to functional data remains unknown.

In a typical functional data setting, a sample of $n$ curves are
observed at a set of discrete points; denote by $m_i$ the number of
observations for curve $i$. The existing literature focuses on two
antithetical data types: the first one, referred to as dense functional
data, is the case where each $m_i$ is larger than some power of~$n$;
the second type, referred to as sparse functional data, is the
situation where each $m_i$ is bounded by a finite positive number or
follows a fixed distribution. The methodologies used to treat the two
situations have been different in the literature. For dense functional
data, the conventional approach is to smooth each individual curve
first before further analysis; see \citet{RS05},
\citet{HMW06} and \citet{ZC07}. For sparse functional data,
limited information is given by the sparsely sampled observations from
each individual curve and hence it is essential to pool the data in
order to conduct inference effectively; see Yao, M\"uller and Wang
(\citeyear{YMW05a}) and \citet{HMW06}. However, in practice it is
possible that some sample curves are densely observed while others are
sparsely observed. Moreover, in dealing with real data, it may even be
difficult to classify which scenario we are faced with and hence to
decide which methodology to use.

This paper is aimed at resolving the issues raised in the previous two
paragraphs. The precise goals will be stated after we introduce the
notation in Section~\ref{sec:model}. In a nutshell, we will consider
uniform rates of convergence of the mean and the covariance functions,
as well as rates in the ensuing FPCA, using local-linear smoothers
[\citet{FG95}]. The rates that we obtain will address all possible
scenarios of the $m_i$'s, and we show that the optimal rates for dense
and sparse functional data can be derived as special cases.

This paper is organized as follows. In Section~\ref{sec:model}, we
introduce the model and data structure as well as all of the estimation
procedures. We describe the asymptotic theory of the procedures in
Section~\ref{sec:asymp}, where we also discuss the results and their
connections to prominent results in the literature. Some simulation
studies are provided in Section~\ref{sec:simulation}, and all proofs
are included in Section~\ref{sec:proof}.

%
\section{Model and methodology}\label{sec:model}

Let $\{X(t), t\in[a,b]\}$ be a stochastic process defined on a fixed
interval $[a,b]$.
Denote the mean and covariance function of the process by
\[
\mu(t)=\mathbb{E}\{X(t)\}, \qquad  R(s,t)=\operatorname{cov}\{
X(s),X(t)\},
\]
which are assumed to exist. Except for smoothness conditions on $\mu$
and $R$, we do not impose any parametric structure on the distribution
of $X$. This is a commonly considered situation in functional data
analysis.

Suppose we observe
\[
Y_{ij}=X_i(T_{ij})+U_{ij}, \qquad  i=1,\ldots,n, j=1,\ldots,m_i,
\]
where the $X_i$'s are independent realizations of $X$, the $T_{ij}$'s
are random observational points with density function $f_T(\cdot)$, and
the $U_{ij}$'s are identically\vspace*{2pt} distributed random errors with mean zero
and finite variance $\sigma^2$. Assume that the $X_i$'s, $T_{ij}$'s and
$U_{ij}$'s are all independent. Assume that $m_i\ge2$ and let
$N_i=m_i(m_i-1)$.

Our approach is based on the local-linear smoother; see, for example,
\citet{FG95}. Let $K(\cdot)$ be a symmetric probability density
function on $[0,1]$ and $K_h(t)=(1/h) K(t/h)$ where $h$ is bandwidth.
A~local-linear estimator of the mean function is given by $\widehat \mu
(t)=\widehat a_0 $, where
\[
(\widehat a_0, \widehat a_1)=\mathop{\arg\min}_{a_0, a_1}
{1\over n} \sum_{i=1}^n {1\over m_i} \sum_{j=1}^{m_i} \{
Y_{ij}-a_0-a_1(T_{ij}-t)\}^2
K_{h_\mu}(T_{ij}-t).
\]
%
It is easy to see that
%
%
\begin{equation}\label{eq:local_linear_soluation}
\widehat\mu(t) = {R_0S_2-R_1S_1\over S_0S_2-S_1^2},
\end{equation}
where
\begin{eqnarray*}
S_r&=&{1\over n}\sum_{i=1}^n {1\over m_i } \sum_{j=1}^{m_i} K_{h_\mu
}(T_{ij}-t)
\{(T_{ij}-t)/h_{\mu}\}^r,\\
R_r&=&{1\over n} \sum_{i=1}^n {1\over m_i} \sum_{j=1}^{m_i}
K_{h_\mu}(T_{ij}-t)
\{(T_{ij}-t)/h_\mu\}^r Y_{ij}.
%
\end{eqnarray*}

To estimate $R(s,t)$, we first estimate $C(s,t):=\mathbb{E}\{X(s)X(t)\}$.
Let $\widehat C(s,t)=\widehat a_0$, where
%
%
\begin{eqnarray}\label{e:llsC}
&&(\widehat a_0,\widehat a_1,\widehat a_2) \nonumber\\
&&\qquad= \mathop{\arg\min}_{a_0,a_1,a_2}
{1\over n} \sum_{i=1}^n  \biggl[ {1\over N_i}
\sum_{k\not=j} \{Y_{ij}Y_{ik}-a_0\nonumber\\[-8pt]\\[-8pt]
&&\qquad\quad\hspace*{96.31pt}{} -a_1(T_{ij}-s)-a_2(T_{ik}-t)\}^2
\nonumber\\
&&\qquad\quad\hspace*{100.39pt}{} \times K_{h_R}(T_{ij}-s) K_{h_R}(T_{ik}-t) \biggr],\nonumber
\end{eqnarray}
with $\sum_{k\not=j}$ denoting sum over all $k,j=1,\ldots,m_i$ such
that $k\not=j$.
It follows that
\[
\widehat C(s,t)=(\mathcal{A}_1 R_{00}-\mathcal{A}_2 R_{10}-\mathcal
{A}_3 R_{01}) \mathcal{B}^{-1},
\]
where
\begin{eqnarray*}
\mathcal{A}_1&=&S_{20}S_{02}-S_{11}^2, \qquad \mathcal{A}_2=S_{10}
S_{02}-S_{01}S_{11}, \qquad \mathcal{A}_3=S_{01}S_{20} -S_{10}S_{11},\\
\mathcal{B}&=&\mathcal{A}_1 S_{00}-\mathcal{A}_2 S_{10}-\mathcal
{A}_3 S_{01},\\
S_{pq}&=&{1\over n} \sum_{i=1}^n {1\over N_i}\sum_{k\not=j}
\biggl({T_{ij}-s
\over h_R} \biggr)^p  \biggl({T_{ik}-t
\over h_R} \biggr)^q K_{h_R}(T_{ij}-s)
K_{h_R}(T_{ik}-t),\\
R_{pq}&=&{1\over n}\sum_{i=1}^n{1\over N_i}
\sum_{k\not=j} Y_{ij}Y_{ik}  \biggl({T_{ij}-s
\over h_R} \biggr)^p  \biggl({T_{ik}-t
\over h_R} \biggr)^q K_{h_R}(T_{ij}-s)
K_{h_R}(T_{ik}-t).
\end{eqnarray*}
We then estimate $R(s,t)$ by
%
%
\begin{equation}\label{eq:cov_est_sol}
\widehat R(s,t) = \widehat C(s,t)-\widehat\mu(s)\widehat\mu(t).\vadjust{\goodbreak}
\end{equation}

To estimate $\sigma^2$, we first estimate
$V(t) := C(t,t)+\sigma^2$ by $\widehat V(t) = \widehat a_0$, where
\[
(\widehat a_0, \widehat a_1)=\mathop{\arg\min}_{a_0, a_1}
{1\over n} \sum_{i=1}^n {1\over m_i} \sum_{j=1}^{m_i} \{
Y_{ij}^2-a_0-a_1(T_{ij}-t)\}^2
K_{h_V}(T_{ij}-t).
\]
As in (\ref{eq:local_linear_soluation}),
%
%
\begin{equation}\label{eq:lls_error-var}
\widehat V(t) = {Q_0S_2-Q_1S_1\over S_0S_2-S_1^2},
\end{equation}
where
\[
Q_r= {1\over n} \sum_{i=1}^n {1\over m_i} \sum_{j=1}^{m_i}
K_{h_V}(T_{ij}-t)
\{(T_{ij}-t)/h_V\}^r Y_{ij}^2.
%
\]
We then estimate $\sigma^2$ by
\[
\widehat\sigma^2 = {1\over b-a} \int_a^b \{\widehat V(t)-\widehat
C(t,t)\}\, dt.
\]

For the problem of mean and covariance estimation, the literature has
focused on dense and sparse functional data. The sparse case roughly
refers to the situation where each $m_i$ is essentially bounded by some
finite number $M$. Yao, M\"uller and Wang (\citeyear{YMW05a}) and \citet{HMW06} considered this case and also
used local-linear smoothers in their estimation procedures. The
difference between the estimators in (\ref{eq:local_linear_soluation}),
(\ref{eq:cov_est_sol}) and those considered in Yao, M\"uller and Wang
(\citeyear{YMW05a}) and \citet{HMW06} is
essentially that we attach weights, $m_i^{-1}$ and $N_i^{-1}$, to each
curve $i$ in the optimizations [although Yao, M\"uller and Wang
(\citeyear{YMW05a}) smoothed the residuals in
estimating $R$]. One of the purposes of those weights is to ensure that
the effect that each curve has on the optimizers is not overly affected
by the denseness of the observations.

Dense functional data roughly refer to data for which each $m_i\ge
M_n\to\infty$ for some sequence $M_n$, where specific assumptions on
the rate of increase of $M_n$ are required for this case to have a
distinguishable asymptotic theory in the estimation of the mean and
covariance. \citet{HMW06} and \citet{ZC07} considered the
so-called ``smooth-first-then-estimate'' approach, namely, the approach
that first preprocesses the discrete functional data by smoothing, and
then adopts the empirical estimators of the mean and covariance based
on the smoothed functional data. See also \citet{RS05}.

As will be seen, our approach is suitable for both sparse and dense
functional data. Thus, one particular advantage is that we do not have
to discern data type---dense, sparse or mixed---and
decide which methodology should be used accordingly.
In Section~\ref{sec:asymp}, we will provide the convergence rates of $\widehat\mu(t),
\widehat R(s,t)$ and $\widehat\sigma^2$,
and also those of the estimated eigenvalues and eigenfunctions of the
covariance operator of $X$. The novelties of our results include:
\begin{enumerate}[(a)]
\item[(a)]
Almost-sure uniform rates of convergence for $\widehat\mu(t)$ and
$\widehat
R(s,t)$ over the entire
range of $s,t$ will be proved.
\item[(b)]
The sample sizes $m_i$ per curve will be completely flexible. For the
special cases of
dense and sparse functional data, these rates match the best
known/conjectured rates.
\end{enumerate}

\section{Asymptotic theory} \label{sec:asymp}

To prove a general asymptotic theory,
assume that $m_{i}$ may depend on $n$ as well, namely, $m_i=m_{in}$.
However, for simplicity we continue to
use the notation $m_i$. Define
\[
\gamma_{nk} =  \Biggl(n^{-1}\sum_{i=1}^nm_i^{-k} \Biggr)^{-1},\qquad
k=1,2,\ldots,
\]
which is the $k$th order harmonic mean of $\{m_i\}$, and for any
bandwidth $h$,
\[
\delta_{n1}(h)=  [\{1+(h\gamma_{n1})^{-1}\}\log n/n ]^{1/2}
\]
and
\[
\delta_{n2}(h)= [\{1+(h\gamma_{n1})^{-1} + (h^2\gamma
_{n2})^{-1}\}\log n/n ]^{1/2}.
\]

We first state the assumptions. In the following $h_\mu, h_R$ and $h_V$
are bandwidths, which are assumed to change with $n$.
\begin{enumerate}[(C1)]
\item[(C1)]
For some constants $m_T>0$ and $M_T<\infty$, $m_T\le f_T(t)\le M_T$
for all $t\in[a,b]$.
Further, $f_T$ is differentiable with a bounded derivative.
\item[(C2)]
The kernel function $K(\cdot)$ is a
symmetric probability
density function on $[-1,1]$, and is of bounded variation on $[-1,1]$.
Denote $\nu_2=\int_{-1}^1 t^2 K(t) \,dt$.
\item[(C3)]
$\mu(\cdot)$ is twice differentiable and the second derivative is
bounded on
$[a,b]$.
\item[(C4)]
All second-order partial derivatives of $R(s,t)$ exist and are bounded
on $[a,b]^2$.
\item[(C5)]
$\mathbb{E}(|U_{ij}|^{\lambda_\mu})<\infty$ and
$\mathbb{E}(\sup_{t\in[a,b]}|X(t)|^{\lambda_\mu})<\infty$
for some $\lambda_\mu\in(2,\infty)$;
$h_\mu\to0$ and $(h_\mu^2 + h_\mu/\gamma_{n1})^{-1}
(\log n/ n)^{1-2/{\lambda_\mu}} \to0$ as $n\to\infty$.
\item[(C6)]
$\mathbb{E}(|U_{ij}|^{2\lambda_R})<\infty$ and
$\mathbb{E}(\sup_{t\in[a,b]}|X(t)|^{2\lambda_R})<\infty$
for some $\lambda_R\in(2,\infty)$;
$h_R\to0$ and $(h_R^4+h_R^3/\gamma_{n1}+h_R^2/\gamma_{n2})^{-1}
(\log n/n)^{1-2/\lambda_R} \to0$ as $n\to\infty$.
\item[(C7)]
$\mathbb{E}(|U_{ij}|^{2\lambda_V})<\infty$ and
$\mathbb{E}(\sup_{t\in[a,b]}|X(t)|^{2\lambda_V})<\infty$
for some $\lambda_V\in(2,\infty)$;
$h_V \to0$ and $(h_V^2 + h_V/\gamma_{n1})^{-1}
(\log n/ n)^{1-2/{\lambda_V}} \to0$ as $n\to\infty$.
\end{enumerate}
The moment condition $\mathbb{E}(\sup_{t\in[a,b]}|X(t)|^{\lambda
})<\infty$
in (C5)--(C7) hold rather generally; in particular, it holds for
Gaussian processes with continuous sample paths [cf. \citet{LS70}]
for all $\lambda>0$. This condition was also adopted by
\citet{HMW06}.
%

\subsection{Convergence rates in mean estimation}\label{sec:mean}
The convergence rate of $\widehat\mu(t)$ is given in the following result.
\begin{theorem}\label{thm:ucr}
Assume that \textup{(C1)--(C3)} and \textup{(C5)} hold.
Then
%
%
\begin{equation}
\label{eq:ucr_mu}
{\sup_{t\in[a,b]}} | \widehat\mu(t)-\mu(t)|=O\bigl(h_\mu^{2}+ \delta
_{n1}(h_\mu)\bigr)  \qquad\mbox{a.s.}
\end{equation}
\end{theorem}

The following corollary addresses the special cases of sparse and
dense functional data. For convenience, we use the notation
$a_n\lesssim b_n$ to mean $a_n=O(b_n)$.
\begin{cor}\label{eq:ucr_mu_cor}
Assume that \textup{(C1)--(C3)} and \textup{(C5)} hold.
\begin{enumerate}[(a)]
\item[(a)]
If $\max_{1\le i\le n} m_i\le M$ for some fixed $M$, then
%
%
\begin{equation}\label{eq:ucr_mu_sparse}
{\sup_{t\in[a,b]}} |\widehat\mu(t)-\mu(t)|=O\bigl(h_\mu^{2}+ \{\log
n/(nh_\mu
)\}^{1/2}\bigr)  \qquad\mbox{a.s.}
\end{equation}
\item[(b)]
If $\min_{1\le i\le n} m_i\ge M_n$ for some sequence $M_n$ where
$M_n^{-1} \lesssim\break h_\mu\lesssim(\log n/n)^{1/4}$
is bounded away from $0$, then
\[
{\sup_{t\in[a,b]}} | \widehat\mu(t)-\mu(t)|=O(\{\log n/n\}^{1/2})
\qquad\mbox{a.s.}
\]
\end{enumerate}
\end{cor}

The proofs of Theorem~\ref{thm:ucr}, as the proofs of other results,
will be given in Section~\ref{sec:proof}. First, we give a few remarks on these results.

\subsubsection*{Discussion}
\begin{enumerate}
\item
On the right-hand side of (\ref{eq:ucr_mu}), $O(h_\mu^2)$ is a bound
for bias while
$\delta_{n1}(h_\mu)$ is a bound for ${\sup_{t\in[a,b]}}|\widehat
\mu(t)-
\mathbb{E}(\widehat\mu(t))|$.
The derivation of the bias is easy to understand and is essentially the
same as in classical
nonparametric regression. The derivation of the second bound is more
involved and
represents our main contribution in this result. To obtain a uniform
bound for
$|\widehat\mu(t)- \mathbb{E}(\widehat\mu(t))|$ over $[a,b]$, we
first obtained a
uniform bound over a finite
grid on $[a,b]$, where the grid grows increasingly dense with $n$, and
then showed that the difference between
the two uniform bounds is asymptotic negligible.
This approach was inspired by \citet{HJS88}, which focused on
nonparametric
regression. One of the main difficulties in our result is that we need
to deal within-curve dependence,
which is not an issue in classical nonparametric regression. Note that
the dependence between $X(t)$ and
$X(t')$ typically becomes stronger as $|t-t'|$ becomes smaller. Thus,
for dense functional data, the within-curve
dependence constitutes an integral component of the overall rate derivation.
%
\item
The sparse functional data setting in (a) of Corollary
\ref{eq:ucr_mu_cor} was considered by Yao, M\"uller and Wang (\citeyear{YMW05a})
and \citet{HMW06}. Actually Yao, M\"uller and Wang (\citeyear{YMW05a}) assumes that the
$m_i$'s are i.i.d. positive random variables
with $\mathbb{E}(m_i) < \infty$,
which implies that $0 < 1/\mathbb{E}(m_i) \le\mathbb{E}(1/m_i)\le1$
by Jensen's
inequality; this corresponds to the case
where $\gamma_{n1}$ is bounded
away from 0 and also leads to (\ref{eq:ucr_mu_sparse}).
The rate in (\ref{eq:ucr_mu_sparse}) is the classical nonparametric
rate for estimating a
univariate function. We will refer to this as a one-dimensional rate.
The one-dimensional rate of $\widehat\mu(t)$ was eluded to in
Yao, M\"uller and Wang (\citeyear{YMW05a})
but was not specifically obtained there.
\item
\citet{HMW06} and \citet{ZC07} address the dense
functional data setting in (b) of
Corollary~\ref{eq:ucr_mu_cor}, where both papers take the approach of first
fitting a smooth curve to $Y_{ij}, 1\le j\le m_i$, for each $i$, and
then estimating
$\mu(t)$ and $R(s,t)$ by the sample mean and covariance functions,
respectively, of
the fitted curves. Two drawbacks are:\looseness=-1
\begin{enumerate}[(a)]
\item[(a)]
Differentiability of the sample curves is required. Thus, for instance,
this approach will not
be suitable for the Brownian motion,
which has continuous but nondifferentiable sample paths.

\item[(b)]
The sample curves that are included in the analysis need to be all
densely observed;
those that do not meet the denseness criterion are dropped even though
they may contain
useful information.
\end{enumerate}\looseness=0
Our approach does not require sample-path differentiability and all of
the data are used
in the analysis. It is interesting to note that (b) of Corollary
\ref{eq:ucr_mu_cor}
shows that root-$n$ rate of convergence for $\widehat\mu$ can be achieved
if the number of
observations per sample curve is at least of the order $(n/\log
n)^{1/4}$ while a similar conclusion
was also reached in \citet{HMW06} for the
smooth-first-then-estimate approach.
%
\item
Our nonparametric estimators $\widehat\mu$, $\widehat R$ and
$\widehat V$ are based
local-linear smoothers, but the methodology
and theory can be easily generalized to higher-order local-polynomial smoothers.
By the equivalent kernel theory for local-polynomial smoothing [\citet{FG95}],
higher-order local-polynomial smoothing is asymptotically equivalent to
higher-order kernel smothing.
Therefore, applying higher-order polynomial smoothing will result in
improved rates for the bias under suitable smoothness
assumptions. The rate for the variance, on the other hand, will remain
the same. In our sparse
setting, if $p$th order local polynomial smoothing is applied under
suitable conditions,
for some positive integer $p$, the uniform convergence rate of
$\widehat\mu
(t)$ will become
\[
{\sup_t} |\widehat\mu(t)-\mu(t)|=O\bigl( h_\mu^{2([p/2]+1)} +\delta
_{n1}(h_\mu) \bigr)  \qquad\mbox{a.s.},
\]
where $[a]$ denotes the integer part of $a$. See
\citet{CK03} and \citet{M96} for support of this claim in different
but related contexts.\vadjust{\goodbreak}

\end{enumerate}

%
\subsection{Convergence rates in covariance estimation}\label{sec:cov_est}
The following results give the convergence rates for $\widehat R(s,t)$ and
$\widehat\sigma^2$.
\begin{theorem}\label{thm:ucr_cov} Assume that \textup{(C1)--(C6)} hold.
Then
%
%
\begin{equation}\label{e:Rrate}
{\sup_{s,t\in[a,b]}} | \widehat R(s,t)-R(s,t)|=O\bigl(h_\mu^2 + \delta
_{n1}(h_\mu) + h_R^2
+ \delta_{n2}(h_R)\bigr)\qquad
\mbox{a.s.}\hspace*{-25pt}
\end{equation}
\end{theorem}
\begin{theorem}\label{thm:error_var}
Assume that \textup{(C1), (C2), (C4), (C6)} and \textup{(C7)} hold. Then
%
%
\begin{equation}\label{e:error_var}
\widehat\sigma^2-\sigma^2 = O\bigl( h_R^2+\delta_{n1}(h_R)+\delta
_{n2}^2(h_R)+h_V^2+\delta_{n1}^2(h_V)\bigr)\qquad
\mbox{a.s.}
\end{equation}
\end{theorem}

We again highlight the cases of sparse and dense functional data.
\begin{cor}\label{cor:ucr_cov_cor}
Assume that \textup{(C1)--(C7)} hold.
\begin{enumerate}[(a)]
\item[(a)]
Suppose that $\max_{1\le i\le n} m_i\le M$ for some fixed $M$. If
$h_R^2 \lesssim h_\mu\lesssim h_R$, then
%
%
\begin{equation}\label{eq:ucr_cov_sparse}
{\sup_{s,t\in[a,b]}}| \widehat R(s,t)-R(s,t)|=O\bigl(h_R^{2}+ \{\log
n/(nh_R^2)\}
^{1/2}\bigr)  \qquad\mbox{a.s.}
\end{equation}
If $h_V+(\log n/n)^{1/3} \lesssim h_R \lesssim h_V^2 n/\log n$, then
\[
\widehat\sigma^2-\sigma^2 = O\bigl(h_R^{2}+ \{\log n/(nh_R)\}^{1/2}\bigr)
\qquad\mbox{a.s.}
\]
\item[(b)]
If $\min_{1\le i\le n} m_i\ge M_n$ for\vspace*{2pt} some sequence $M_n$ where
$M_n^{-1}\lesssim h_\mu, h_R, h_V
\lesssim(\log n/n)^{1/4}$, then both ${\sup_{s,t\in[a,b]}} |\widehat
R(s,t)-R(s,t)|$ and $\widehat\sigma^2-\sigma^2$ are\break
$O(\{\log n/n\}^{1/2})$ a.s.
%
%
\end{enumerate}
\end{cor}

\subsubsection*{Discussion}
\begin{enumerate}
\item
The rate in (\ref{eq:ucr_cov_sparse}) is the classical nonparametric
rate for estimating
a surface (bivariate function), which will be referred to as a
two-dimensional rate.
Note $\widehat\sigma^2$ has a one-dimensional rate in the sparse
setting, while
both $\widehat R(s,t)$ and $\widehat\sigma^2$ have root-$n$ rates
in the dense setting.
Most of the discussions in Section~\ref{sec:mean} obviously also apply
here and will not be repeated.
\item
Yao, M\"uller and Wang (\citeyear{YMW05a}) smoothed the products of residuals instead of
$Y_{ij}Y_{ik}$ in the local linear smoothing
algorithm in (\ref{e:llsC}). There is some evidence that a slightly
better rate can be achieved in that procedure.
However, we were not successful in establishing such a rate rigorously.
\end{enumerate}

%
\subsection{Convergence rates in FPCA}\label{sec:fpca_rate}

By (C5), the covariance function has the spectral decomposition
\[
R(s,t)=\sum_{j=1}^\infty\omega_j \psi_j(s) \psi_j(t),
\]
where $\omega_1\ge\omega_2\ge\cdots\ge0$ are the eigenvalues of
$R(\cdot,\cdot)$ and
the $\psi_j$'s are the corresponding eigenfunctions. The $\psi_j$'s
are also known as the functional principal components. Below, we
assume that the nonzero $\omega_j$'s are distinct.\vadjust{\goodbreak}

Suppose $\widehat R(s,t)$ is the covariance estimator given in Section
\ref{sec:model}, and it admits the following spectral decomposition:
\[
\widehat R(s,t)=\sum_{j=1}^\infty\widehat\omega_j \widehat\psi
_j(s) \widehat\psi_j(t),
\]
where $\widehat\omega_1> \widehat\omega_2> \cdots$ are the estimated
eigenvalues and the $\widehat\psi_j$'s are the corresponding estimated
principal components. Computing the eigenvalues and eigenfunctions of
an integral operator with a symmetric kernel is a well-studied problem
in applied mathematics. We will not get into that aspect of FPCA in
this paper.

Notice also that $\psi_j(t)$ and $\widehat\psi_j(t)$ are
identifiable up
to a sign change. As pointed out in \citet{HMW06}, this causes no
problem in practice, except when we discuss the convergence rate of
$\widehat\psi_j$. Following the same convention as in \citet{HMW06},
we let $\psi_j$ take an arbitrary sign but choose $\widehat\psi_j$ such
that $\|\widehat\psi_j-\psi_j\|$ is minimized over the two
signs,\vspace*{1pt}
where $\|f\|:=\{\int f^2(t) \,dt \}^{1/2}$ denotes the usual $L^2$-norm
of a function $f\in L^2[a,b]$.

Below let $j_0$ be a arbitrary fixed positive constant.
\begin{theorem}\label{thm:pca_rate}
Under conditions \textup{(C1)--(C6)}, for $1\le j \le j_0$:
\begin{enumerate}[(a)]
\item[(a)]
$\widehat\omega_j-\omega_j=O((\log n/n)^{1/2}+h_\mu^2+h_R^2+\delta
_{n1}^2(h_\mu)
+\delta_{n2}^2(h_R))$ a.s.;
\item[(b)]
$\|\widehat\psi_j-\psi_j\|=O(h_{\mu}^2+\delta_{n1}(h_\mu)+h_R^2
+\delta_{n1}(h_R)+\delta_{n2}^2(h_R))$ a.s.;
\item[(c)]
${\sup_t }|\widehat\psi_j(t)-\psi_j(t)|=O(h_{\mu}^2+\delta
_{n1}(h_\mu)+h_R^2
+\delta_{n1}(h_R)+\delta_{n2}^2(h_R))$ a.s.
\end{enumerate}
\end{theorem}

Theorem~\ref{thm:pca_rate} is proved by using the asymptotic
expansions of eigenvalues
and eigenfunctions of an estimated covariance function developed by
\citet{HH06}, and by applying the strong uniform
convergence rate of $\widehat R(s,t)$ in Theorem~\ref{thm:ucr_cov}.
In the
special case of sparse and dense functional data, we have the following
corollary.
\begin{cor}\label{cor:pca_rate}
Assume that \textup{(C1)--(C6)} hold. Suppose that\break $\max_{1\le i \le n} m_i\le
M$ for some fixed $M$.
Then the following hold for all $1\le j \le j_0$:
\begin{enumerate}[(a)]
\item[(a)]
If $(\log n/n)^{1/2} \lesssim h_\mu, h_R \lesssim(\log n/n)^{1/4}$ then
$\widehat\omega_j-\omega_j=O(\{\log n/n\}^{1/2})$ a.s.
\item[(b)]
If $h_\mu+(\log n/n)^{1/3} \lesssim h_R \lesssim h_\mu$ then both of
$\|\widehat\psi_j-\psi_j\|$ and ${\sup_t} |\widehat\psi_j(t)-\psi_j(t)|$
have the rate
$O(h_R^{2}+ \{\log n/(nh_R)\}^{1/2})$ a.s.
\end{enumerate}
If $\min_{1\le i\le n} m_i\ge M_n$ for some sequence $M_n$ where
$M_n^{-1}\lesssim h_\mu, h_R
\lesssim(\log n/n)^{1/4}$, then, for $1\le j \le j_0$, all of
$\widehat\omega_j-\omega_j$, $\|\widehat\psi_j-\psi_j\|$ and
$\sup_t |\widehat
\psi_j(t)-\psi_j(t)|$ have the rate
$O(\{\log n/n\}^{1/2})$.
\end{cor}

\subsubsection*{Discussion}
\begin{enumerate}
\item
Yao, M\"uller and Wang (\citeyear{YMW05a}, \citeyear{YMW05b})
developed rate estimates for the quantities in Theorem
\ref{thm:pca_rate}. However, they are not optimal.\vspace*{1pt}
\citet{HMW06} considered the rates of $\widehat\omega_j-\omega
_j$ and
$\|\widehat\psi_j-\psi_j\|$. The most striking insight of their results
is that for sparse functional
data, even though the estimated covariance operator has the two-dimensional
nonparametric rate, $\widehat\psi_j$~converges at a one-dimensional rate
while $\widehat\omega_j$ converges
at a root-$n$ rate if suitable smoothing parameters are used;
remarkably they also established the asymptotic distribution of $\|
\widehat
\psi_j-\psi_j\|$.
At first sight, it may seem counter-intuitive that the convergence
rates of $\widehat\omega_j$
and $\widehat\psi_j$ are faster than that of $\widehat R$, since
$\widehat\omega_j$
and $\widehat\psi_j$ are computed from~$\widehat R$. However, this
can be
easily explained.
For example, by (4.9) of \citet{HMW06}, $\widehat\omega
_j-\omega_j
= \iint(\widehat R(s,t)-R(s,t)) \psi_j(s) \psi_j(t) \,ds \,dt +
\mathrm{lower}$-order terms;
integrating $\widehat R(s,t)-R(s,t)$ in this expression results in
extra smoothing, which leads to a faster convergence rate.
\item
Our almost-sure convergence rates are new. However,
for both dense and sparse functional data, the rates on
$\widehat\omega_j-\omega_j$ and $\|\widehat\psi_j-\psi_j\|$ are slightly
slower than the in-probability convergence
rates obtained in \citet{HMW06}, which do not contain the $\log n$
factor at various places of our rate bounds. This is due to the fact
that our
proofs are tailored to strong uniform convergence rate derivation. However,
the general strategy in our proofs is amenable to deriving
in-probability convergence
rates that are comparable to those in \citet{HMW06}.\looseness=1
\item
A potential estimator the covariance function $R(s,t)$ is
\[
\widetilde R(s,t) := \sum_{j=1}^{J_n} \widehat\omega_j \widehat
\psi_j(s)\widehat\psi_j(t)
\]\looseness=0
for some $J_n$. For the sparse case, in view of the one-dimensional
uniform rate of $\widehat\psi_j(t)$ and
the root-$n$ rates of $\widehat\omega_j$, it might be possible to choose
$J_n\to\infty$ so that $\widetilde R(s,t)$ has a faster
rate of convergence than does $\widehat R(s,t)$. However, that
requires the
rates of $\widehat\omega_j$
and $\widehat\psi_j(t)$ for an unbounded number of $j$'s, which we
do not
have at this point.

\end{enumerate}

The proof of the theorems will be given in Section~\ref{sec:proof}, whereas the
proofs of the corollaries are straightforward and are omitted.

\section{Simulation studies}\label{sec:simulation}
\subsection{Simulation 1}
To illustrate the finite sample performance of the method, we perform a
simulation study. The data are generated from the following model:
\[
Y_{ij}=X_i(T_{ij})+U_{ij}\qquad   \mbox{with }
X_{i}(t)= \mu(t)+ \sum_{k=1}^3 \xi_{ik} \psi_j(t),
\]
where $T_{ij}\sim \operatorname{Uniform}[0,1]$, $\xi_{ik}\sim\operatorname
{Normal}(0,\omega
_j)$ and $U_{ij}\sim\operatorname{Normal}(0,\sigma^2)$ are
independent variables.
Let
\begin{eqnarray*}
&&\mu(t)=5(t-0.6)^2, \qquad \psi_1(t)=1,\\
&& \psi_2(t)=\sqrt{2}\sin(2\pi t), \qquad \psi_3(t)=\sqrt{2}\cos
(2\pi t)
\end{eqnarray*}
and $(\omega_1,\omega_2,\omega_3,\sigma^2)=(0.6,0.3,0.1, 0.2)$.

We let $n=200$ and $m_i=m$ for all $i$. In each simulation run, we
generated 200 trajectories from the model above,
and then we compared the estimation results for $m=5$, $10$, $50$ and
$\infty$.
When $m=\infty$, we assumed that we know the whole trajectory and so
no measurement error was included.
Note that the cases of $m=5$ and $m=\infty$ may be viewed as
representing sparse and complete functional data, respectively,
whereas those of $m=10$ and $m=50$ represent scenarios between the two
extremes. For each $m$ value, we estimated the mean
and covariance functions and used the estimated covariance function to
conduct FPCA. The simulation was then repeated 200 times.\looseness=-1

For $m=5,10,50$, the estimation was carried out as described in
Section~\ref{sec:model}.
For $m=\infty$, the estimation procedure was different since no kernel
smoothing is needed;
in this case, we simply discretized each curve on a dense grid, then
the mean and covariance functions were estimated
using the gridded data.

Notice that $m=\infty$ is the ideal situation where we have the
complete information of each curve, and the estimation results under
this scenario represent the best we can do and all of the estimators
have root-$n$ rates. Our asymptotic theory shows that $m\to\infty$ as
a function of $n$, and if $m$ increases with a fast enough rate, the
convergence rates for the estimators are also root-$n$.
We intend to demonstrate this based on simulated data.

The performance of the estimators depends on the choice of bandwidths
for $\mu(t)$, $C(s,t)$ and $V(t)$, and the best bandwidths
vary with $m$. The bandwidth selection problem turns out to be very
challenging. We have not
come across a data-driven procedure that
works satisfactorily and so this is an important problem for future research.
For lack of a better approach, we tried picking the bandwidths by
the integrated mean square error (IMSE); that is, for each $m$ and for
each function above,
we calculated the IMSE over
a range of $h$ and selected the one that minimizes the IMSE. The
bandwidths picked that way worked quite well for the
inference of the mean, covariance and the leading principal components,
but less well for $\sigma^2$ and the eigenvalues.
After experimenting with a number of bandwidths, we decided to used
bandwidths that are slightly smaller than the ones
picked by IMSE. They are reported in Table~\ref{tb:simu_bandwidth}.
Note that undersmoothing in functional principal component analysis was
also advocated by \citet{HMW06}.\vadjust{\goodbreak}

%
%
\begin{table}
\tablewidth=250pt
\caption{Bandwidths in simulation 1}\label{tb:simu_bandwidth}
\begin{tabular*}{\tablewidth}{@{\extracolsep{\fill}}lccc@{}}
\hline
& $\bolds{h_\mu}$ & $\bolds{h_R}$ & $\bolds{h_V}$ \\
\hline
$m=5$ & 0.153 & 0.116 & 0.138\\
$m=10$ & 0.138 & 0.103 & 0.107\\
$m=50$ & 0.107 & 0.077 & 0.084\\
\hline
\end{tabular*}
\end{table}
%

%
\begin{figure}[b]

\includegraphics{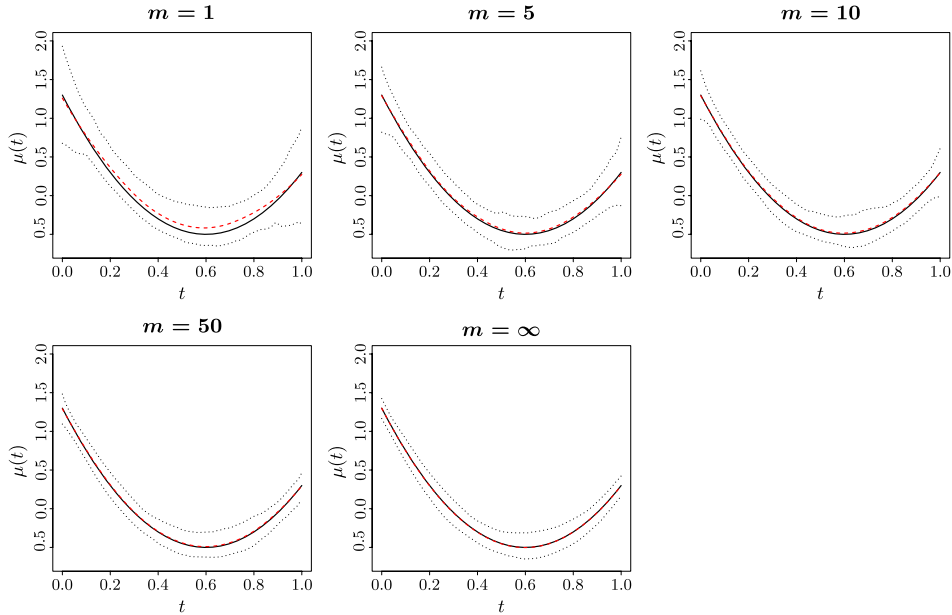}

\caption{Estimated mean function in simulation 1. In each panel, the
solid line is the true mean
function, the dashed line is the pointwise mean and the two dotted
lines are the pointwise
$1\%$ and $99\%$ percentiles of the estimator of the mean function
based on 200 runs.}
\label{fig:simu_mean}
\end{figure}

\begin{figure}

\includegraphics{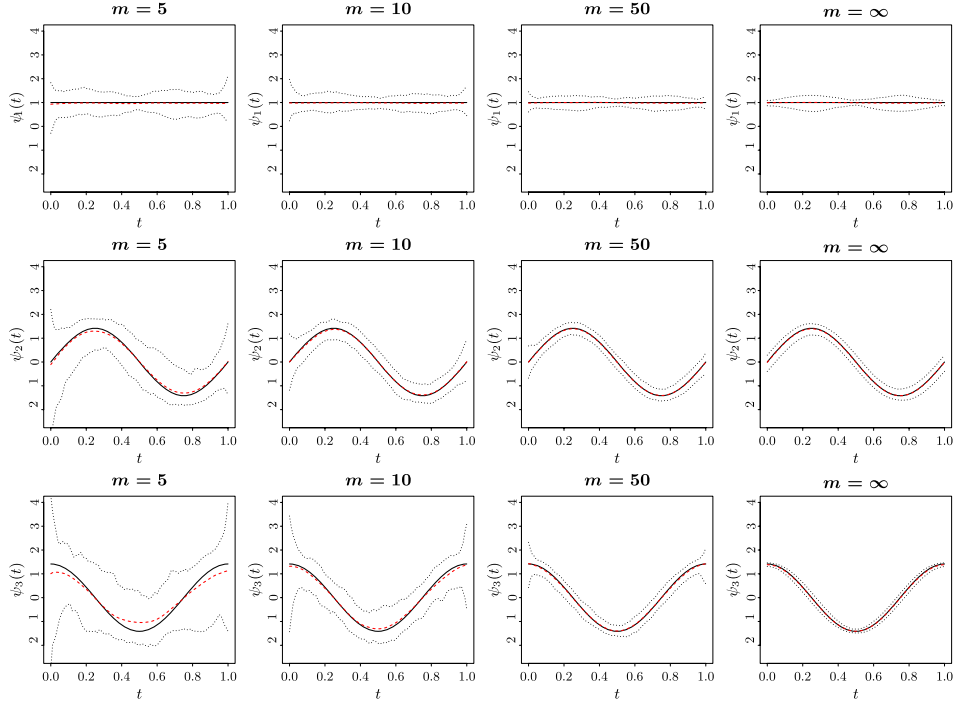}

\caption{Estimated eigenfunctions in simulation 1. In
each panel, the solid line is the eigenfunction, the dashed line is the
pointwise mean and the two dotted lines are the pointwise $1\%$ and
$99\%$ percentiles of the estimator of the eigenfunction in 200 runs.
The three rows correspond to $\psi_1$, $\psi_2$ and $\psi_3$; different
columns correspond to different $m$ values.}
\label{fig:simu_eigenfunction}
\end{figure}

\begin{figure}

\includegraphics{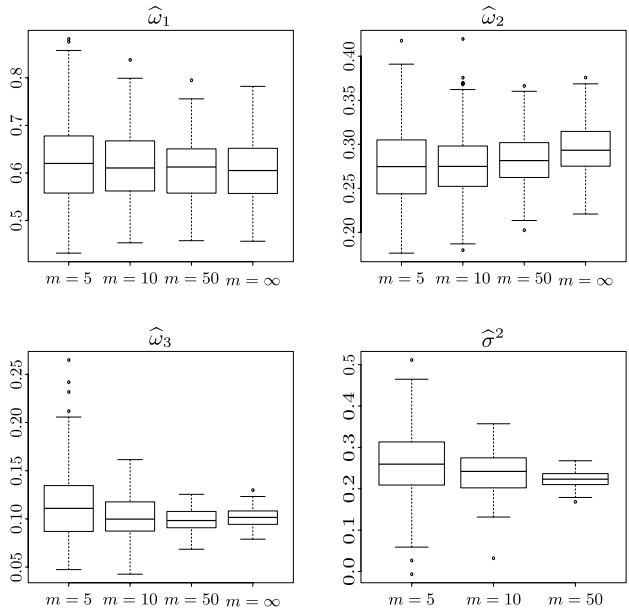}

\caption{Box plots for $\widehat\omega_1$, $\widehat\omega_2$,
$\widehat\omega_3$ and $\widehat\sigma^2$ in simulation 1.}
\label{fig:simu_eigenvalue_var}\vspace*{-3pt}
\end{figure}

The estimation results for $\mu(\cdot)$ are summarized in Figure
\ref{fig:simu_mean},
where we plot the mean and the pointwise first and 99th
percentiles of the estimator.
To compare with standard nonparametric regression, we also provide the
estimation results for $\mu$ when $m=1$; note that in this case the
covariance function is not estimable since there is
no within-curve information. As can be seen, the estimation result for
$m=1$ is not very different from
that of $m=5$, reconfirming the nonparametric convergence rate of
$\widehat
\mu$ for sparse functional data.
It is somewhat difficult to describe the estimation results of the
covariance function
directly. Instead, we summarize the results on $\psi_k(\cdot)$ and
$\omega_k$
in Figure~\ref{fig:simu_eigenfunction}, where we plot the mean
and the pointwise first and 99th percentiles
of the estimated eigenfunctions.
In Figure~\ref{fig:simu_eigenvalue_var}, we also show the empirical
distributions of
$\widehat\omega_k$ and $\widehat\sigma^2$.
In all of the scenarios, the performance of the estimators improve with
$m$; by $m=50$, all of the the estimators perform
almost as well as those for $m=\infty$.

%

%

\subsection{Simulation 2}
To illustrate that the proposed methods are applicable even to the
cases that the trajectory of $X$ is not smooth, we now present a second
simulation study where $X$ is standard Brownian motion. Again, we set
the time window $[a,b]$ to be $[0,1]$. It is well known that the
covariance function of $X$ is
$R(s,t)=\min(s,t), s,t\in[0,1]$, which has an infinite spectral
decomposition with
\[
\omega_k=4/ \{(2k-1)^2 \pi^2\}, \qquad \psi_k(t)=\sqrt{2} \sin\{
(k-1/2) \pi t\}, \qquad  k=1,2,\ldots.
\]
Again, let the observation times be $T_{ij}\sim\operatorname{Uniform}[0,1]$,
$Y_{ij}=X_i(T_{ij})+U_{ij}$, $U_{ij}\sim\operatorname
{Normal}(0,\sigma^2)$. We let
$\sigma^2=0.1^2$, which is comparable to $\omega_3$.

Since $X$ is not differentiable with probability one, smoothing
individual trajectories is not sensible even for large $m$ values.
Also, $R(s,t)$ is not differentiable on the diagonal $\{s=t\}$, and
therefore the smoothness assumption in our theory is not satisfied.
Nevertheless, as we will show below, the proposed method still works
reasonably well. The reason is that the smoothness assumption on
$R(s,t)$ in our theory is meant to guarantee the best convergence rate
for the $\widehat R(s,t)$. When the assumption is mildly violated, the
estimator may still perform well
overall but may have a slower convergence rate at the nonsmooth points.
A similar phenomenon was observed in \citet{Letal07}, which studied
kernel estimation of a stationary covariance function in a time-series setting.

\begin{figure}

\includegraphics{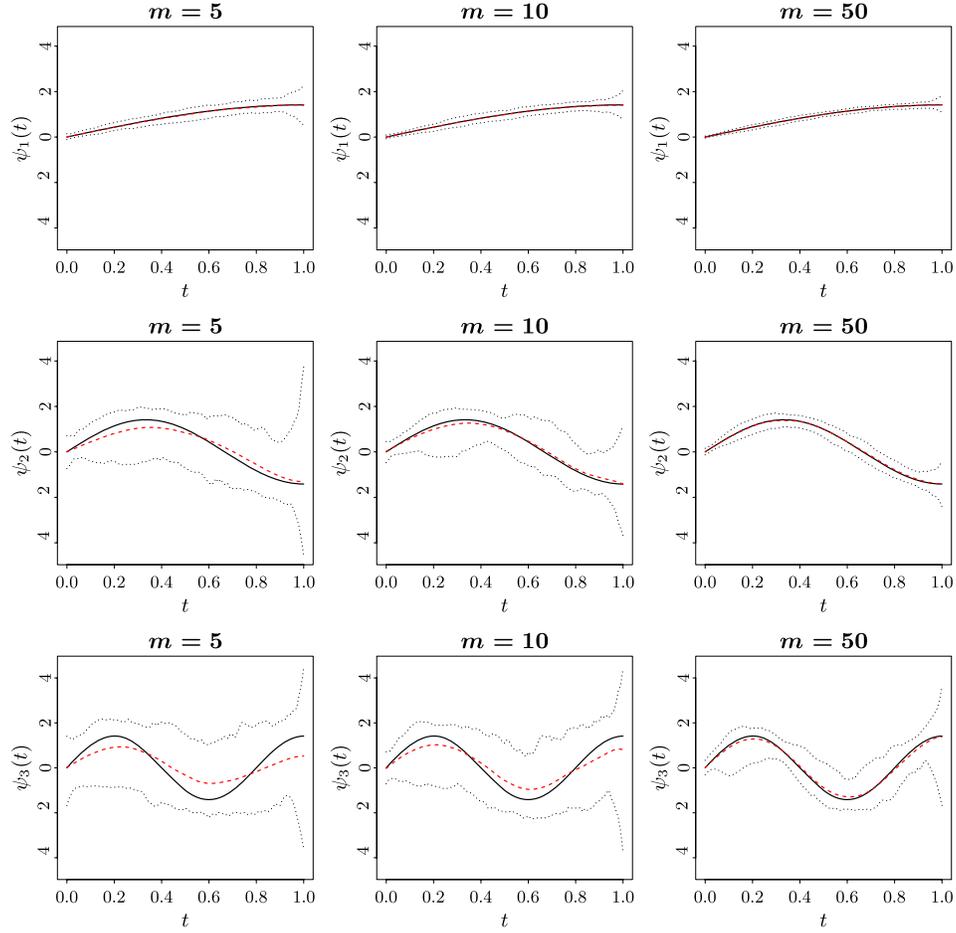}

\caption{Estimated eigenfunctions in simulation 2. In
each panel, the solid line is the eigenfunction, the dashed line is the
pointwise mean and the two dotted lines are the pointwise $1\%$ and
$99\%$ percentiles of the estimator of the eigenfunction in 200 runs.
The three rows correspond to $\psi_1$, $\psi_2$ and $\psi_3$; different
columns correspond to different $m$ values.}
\label{fig:simu2_eigenfunction}
\end{figure}

\begin{figure}

\includegraphics{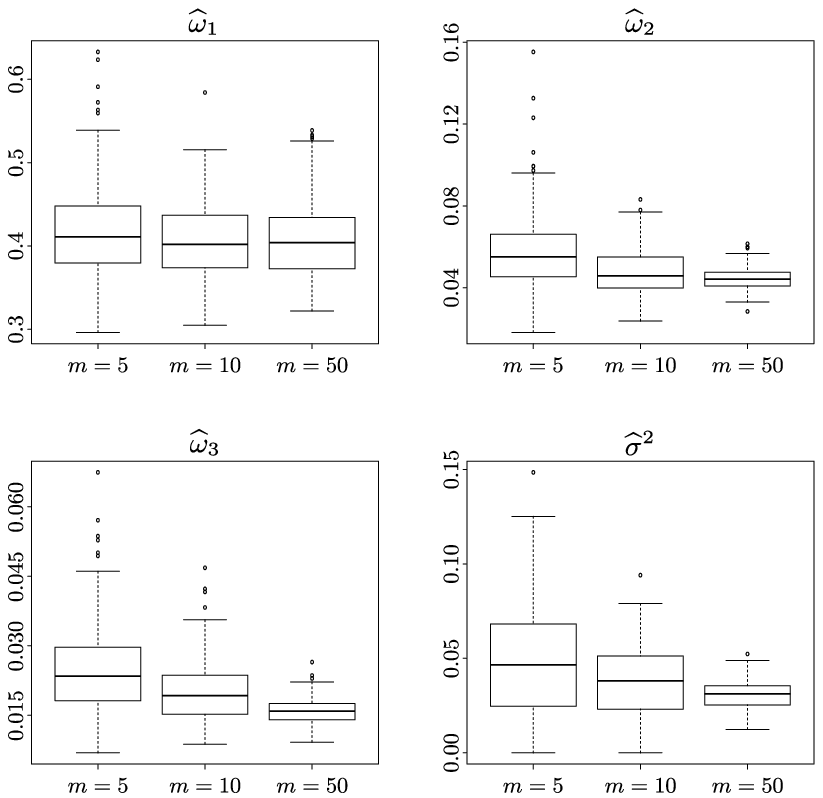}

\caption{Box plots for $\widehat\omega_1$, $\widehat\omega_2$,
$\widehat\omega_3$ and $\widehat\sigma^2$ in simulation 2.}
\label{fig:simu2_eigenvalue_var}
\end{figure}

We set $n=200$ and $m=5$, $10$ or $50$ in our simulations. The
estimation results for the first three eigenfunctions are presented in
Figure~\ref{fig:simu2_eigenfunction}. Again, we plot the mean and the
pointwise first and 99th percentiles
of the estimated eigenfunctions. As can be seen, it is in general much
harder to estimate the higher-order eigenfunctions, and the results
improve as we increase $m$. The empirical distribution of the estimated
eigenvalues as well as $\widehat\sigma^2$ are summarized in Figure
\ref{fig:simu2_eigenvalue_var}. The estimated eigenvalues should be
compared with the true ones, which are $(0.405, 0.045, 0.016)$. When
$m$ is large, the estimated eigenvalues are very close to the true values.

%

%

%
\section{Proofs}\label{sec:proof}
%

\subsection{\texorpdfstring{Proof of Theorem \protect\ref{thm:ucr}}{Proof of Theorem 3.1}}

The proof is an adaptation of familiar lines of proofs established in
nonparametric function literature;
see Claeskens and Van Keilegon (\citeyear{CK03})\vadjust{\goodbreak} and \citet{HJS88}.
For simplicity, throughout this subsection, we abbreviate $h_\mu$ as $h$.
Below, let $t_1\wedge t_2=\min(t_1,t_2)$ and $t_1\vee t_2 = \max(t_1,t_2)$.
Also define $K_{(\ell)}(t)=t^\ell K(t)$ and $K_{h,(\ell)}(v) =
(1/h)K_{(\ell)}(v/h)$.
\begin{lm}\label{lm:lemma1}
Assume that
%
%
\begin{equation}\label{e:lmaassumption}\qquad
\mathbb{E} \Bigl({\sup_{t\in[a,b]} }|X(t)|^{\lambda} \Bigr) <\infty
\quad\mbox{and}\quad \mathbb{E}|U|^{\lambda} < \infty
\qquad\mbox{for some $\lambda\in(2,\infty)$}.
\end{equation}
Let\vspace*{1pt} $\mathcal{Z}_{ij} = X_i(T_{ij})$ or $U_{ij}$ for $1\le i\le n,
1\le j\le m_i$.
Let $c_n$ be any positive sequence tending to $0$ and
$\beta_{n} = c_n^2 + c_n/\gamma_{n1}$.
Assume that $\beta_{n}^{-1}(\log n/n)^{1-2/\lambda} = o(1)$.
Let
%
%
\begin{eqnarray}\label{e:Gn1}
G_n(t_1,t_2) &=& {1\over n}\sum_{i=1}^n  \Biggl\{{1\over m_i} \sum
_{j=1}^{m_i}\mathcal{Z}_{ij} I (T_{ij}\in
[t_1\wedge t_2,t_1\vee t_2] ) \Biggr\}, \nonumber\\[-8pt]\\[-8pt]
G(t_1,t_2)&=&\mathbb{E}\{G_n(t_1,t_2)\} \nonumber
\end{eqnarray}
and
\[
V_n(t,c)={\sup_{|u|\le c}} |G_n(t,t+u)-G(t,t+u)|,\qquad c >0.
\]
Then
%
%
\begin{equation}\label{eq:ucr_G}
\sup_{t\in[a,b]} V_n(t,c_n)=O(n^{-1/2} \{\beta_n \log n\}
^{1/2}) \qquad\mbox{a.s.}
\end{equation}
\end{lm}
\begin{pf}
We can obviously treat the positive and negative parts of $\mathcal{Z}_{ij}$
separately, and will
assume below that $\mathcal{Z}_{ij}$ is nonnegative.
Define an equally-spaced grid $\mathcal{G}:=\{v_k\}$, with $v_k=a+k
c_n$, for
$k=0,\ldots
,[(b-a)/c_n]$, and \mbox{$v_{[(b-a)/c_n]+1}=b$}, where $[\cdot]$ denotes
the greatest integer part. For any $t\in[a,b]$ and $|u|\le c_n$, let
$v_k$ be a grid point that
is within $c_n$ of both $t$ and $t+u$, which exists.
Since
\begin{eqnarray*}
|G_n(t,t+u)-G(t,t+u)| &\le& |G_n(v_k,t+u)-G(v_k,t+u)|\\
&&{} +|G_n(v_k,t)-G(v_k,t)|,
\end{eqnarray*}
we have
\[
|G_n(t,t+u)-G(t,t+u)|\le2 \sup_{t\in\mathcal{G}}V_n(t,c_n).
\]
Thus,
%
%
\begin{equation}\label{eq:a2}
\sup_{t\in[a,b]} V_n(t,c_n)\le2 \sup_{t\in\mathcal{G}} V_n(t,c_n).
\end{equation}
From now on, we focus on the right-hand side of (\ref{eq:a2}).
Let
%
%
\begin{equation}\label{e:qncnan2}
a_n=n^{-1/2} \{\beta_n \log n\}^{1/2}  \quad\mbox{and}\quad  Q_n =
\beta_n/a_n,
\end{equation}
and define $G_n^\ast(t_1,t_2),G^\ast(t_1,t_2)$ and $V_n^\ast(t,c_n)$
in the same way as $G_n(t_1,t_2)$, $G(t_1,t_2)$ and $V_n(t,c_n)$, respectively,
except with $\mathcal{Z}_{ij}I(\mathcal{Z}_{ij}\le Q_n)$ replacing
$\mathcal{Z}_{ij}$.
Then
%
%
\begin{equation}\label{eq:a3}
\sup_{t\in\mathcal{G}} V_n(t,c_n)\le\sup_{t\in\mathcal{G}}
V_n^\ast(t,c_n)+ A_{n1}+A_{n2},
\end{equation}
where
\begin{eqnarray*}
A_{n1}&=&\sup_{t\in\mathcal{G}} \sup_{|u|\le c_n}
\bigl(G_n(t,t+u)-G_n^\ast(t,t+u)\bigr),\\
A_{n2}&=&\sup_{t\in\mathcal{G}} \sup_{|u|\le c_n} \bigl(G(t,t+u)-G^\ast(t,t+u)\bigr).
\end{eqnarray*}
We first consider $A_{n1}$ and $A_{n2}$.
It follows that
%
%
\begin{equation}\label{e:anQn}
a_n^{-1}Q_n^{1-\lambda} = \{\beta_n^{-1}(\log n/n)^{1-2/\lambda}\}
^{\lambda/2} = o(1).
\end{equation}
For all $t$ and $u$, by Markov's inequality,
\begin{eqnarray*}
&&a_n^{-1}\bigl(G_n(t,t+u)-G_n^\ast(t,t+u)\bigr) \\
&&\qquad\le a_n^{-1}{1\over n}\sum_{i=1}^n  \Biggl\{{1\over m_i}
\sum_{j=1}^{m_i} \mathcal{Z}_{ij}I(\mathcal{Z}_{ij}>Q_n) \Biggr\}
\\
&&\qquad\le a_n^{-1} Q_n^{1-\lambda}{1\over n}\sum_{i=1}^n  \Biggl\{{1\over m_i}
\sum_{j=1}^{m_i} \mathcal{Z}_{ij}^{\lambda} I(\mathcal
{Z}_{ij}>Q_n) \Biggr\}
\\
&&\qquad\le a_n^{-1} Q_n^{1-\lambda}{1\over n}\sum_{i=1}^n  \Biggl\{{1\over m_i}
\sum_{j=1}^{m_i} \mathcal{Z}_{ij}^{\lambda}  \Biggr\}.
\end{eqnarray*}
Consider the case $\mathcal{Z}_{ij}=X_i(T_{ij})$, the other case being simpler.
It follows that
\[
{1\over m_i} \sum_{j=1}^{m_i} \mathcal{Z}_{ij}^{\lambda} \le W_i
\qquad \mbox{where }
W_i = {\sup_{t\in[a,b]}} |X_i(t)|^{\lambda}.
\]
Thus,
%
%
\begin{equation}\label{e:GG*}
a_n^{-1}\bigl(G_n(t,t+u)-G_n^\ast(t,t+u)\bigr)
\le a_n^{-1} Q_n^{1-\lambda}{1\over n}\sum_{i=1}^n W_i.
%
%
\end{equation}
By the SLLN, $n^{-1} \sum_{i=1}^n W_i\stackrel{\mathrm
{a.s.}}{\longrightarrow}\mathbb{E}(\sup_{t\in[a,b]}
|X(t)|^{\lambda}) <\infty$.
By (\ref{e:anQn}) and (\ref{e:GG*}), $a_n^{-1} A_{n1}\stackrel{\mathrm
{a.s.}}{\longrightarrow}0$. By (\ref
{e:anQn}) and (\ref{e:GG*}) again,
$a_n^{-1}A_{n2}=0$, and so we have proved
%
%
\begin{equation}\label{eq:a4}
\lim_{n\to\infty} (A_{n1}+A_{n2}) = o(a_n)  \qquad\mbox{a.s.}
\end{equation}
To bound $V_n^\ast(t, c_n)$ for a fixed $t\in\mathcal{G}$, we perform a
further partition. Define $w_n=[Q_n c_n/a_n+1]$ and $u_r=r c_n/w_n$,
for $r=-w_n,-w_n+1,\ldots,w_n$. Note that $G_n^\ast(t,t+u)$ is monotone
in $|u|$ since $\mathcal {Z}_{ij}\ge0$. Suppose that $0\le u_r\le u\le
u_{r+1}$. Then
\begin{eqnarray*}
&& G_n^\ast(t,t+u_r) - G^\ast(t,t+u_r)
+ G^\ast(t,t+u_r) - G^\ast(t,t+u_{r+1}) \\
&&\qquad\le G_n^\ast(t,t+u)-G^\ast(t,t+u) \\
&&\qquad\le G_n^\ast(t,t+u_{r+1}) - G^\ast(t,t+u_{r+1})
+ G^\ast(t,t+u_{r+1}) - G^\ast(t,t+u_r),
\end{eqnarray*}
from which we conclude that
\[
|G_n^\ast(t,t+u)-G^\ast(t,t+u)| \le\max(\xi_{nr},\xi_{n,r+1}) +
G^\ast(t+u_r,t+u_{r+1}),
\]
where
\[
\xi_{nr}=|G_n^\ast(t,t+u_r)-G^\ast(t,t+u_r)|.
\]
The same holds if $u_r\le u\le u_{r+1}\le0$.
Thus,
\[
V_n^\ast(t, c_n)\le\max_{-w_n\le r\le w_n} \xi_{nr} +\max_{-w_n\le
r\le
w_n} G^\ast(t+u_r,t+u_{r+1}).
\]
For all $r$,
\begin{eqnarray*}
G^\ast(t+u_r,t+u_{r+1})
&\le& Q_n \mathbb{P}(t+u_r\le T\le t+u_{r+1})\\
&\le& M_T Q_n (u_{r+1}-u_r) \le M_T a_n.
\end{eqnarray*}
%
Therefore, for any $B$,
%
%
\begin{equation}\label{e:vn*}
\mathbb{P}\{V_n^\ast(t,c_n)\ge B a_n\} \le\mathbb{P}\Bigl\{\max
_{-w_n\le r\le w_n}
\xi_{nr}\ge(B-M_T) a_n\Bigr\}.
%
\end{equation}
Now let $Z_{i}= m_i^{-1} \sum_{j=1}^{m_i} \mathcal{Z}_{ij} I(\mathcal
{Z}_{ij}\le
Q_n) I(T_{ij}\in
(t,t+u_r])$ so that $\xi_{nr}=|{1\over n} \times \sum_{i=1}^n
\{Z_i-\mathbb{E}(Z_i)\}|$. We have $|Z_i-\mathbb{E}(Z_i)|\le Q_n$, and
\[
\sum_{i=1}^n\operatorname{var}(Z_i)\le\sum_{i=1}^n\mathbb
{E}Z_i^2\le M \sum_{i=1}^n
(c_n^2+c_n/m_i)
\le M n \beta_n
\]
for some finite $M$. By Bernstein's inequality,
\begin{eqnarray*}
\mathbb{P}\{\xi_{nr}\ge(B-M_T)a_n\}&\le& \exp \biggl\{ - {(B-M_T)^2
n^2 a_n^2
\over2\sum_{i=1}^n \operatorname{var}(Z_i)+ (2/3)(B-M_T)Q_n na_n}
\biggr\}\\
&\le& \exp \biggl\{- {(B-M_T)^2 n^2 a_n^2
\over2 M n \beta_n+ (2/3)(B-M_T)n
\beta_n} \biggr\} \le n^{-B^\ast},
\end{eqnarray*}
where $B^\ast={(B-M_T)^2 \over2 M +(2/3)(B-M_T)}$.
By (\ref{e:vn*}) and Boole's inequality,
\[
\mathbb{P} \Bigl\{\sup_{t\in\mathcal{G}} V_n^\ast(t,c_n)\ge B
a_n \Bigr\} \!\le\!
\biggl( \biggl[{b-a
\over c_n} \biggr]+1 \biggr)\!  \biggl(2 \biggl[{Q_nc_n \over a_n}+1
\biggr]+1 \biggr) n^{-B^\ast}
\!\le\! C {Q_n\over a_n} n^{-B^\ast}
\]
for some finite $C$. Now $Q_n/a_n = \beta_n/a_n^2 = n/\log n$.
%
So $\mathbb{P}\{V_n^\ast(t,c_n)\ge B a_n\}$ is summable in $n$ if we select
$B$ large enough
such that $B^\ast>2$. By the Borel--Cantelli lemma,
%
%
\begin{equation}\label{eq:a7}
\sup_{t\in\mathcal{G}} V_n^\ast(t,c_n)=O(a_n)  \qquad\mbox{a.s.}
\end{equation}
Hence, (\ref{eq:ucr_G}) follows from combining (\ref{eq:a2}), (\ref
{eq:a3}), (\ref{eq:a4}) and (\ref{eq:a7}).
\end{pf}
\begin{lm}\label{lm:urc_general_kernel}
Let $\mathcal{Z}_{ij}$ be as in Lemma~\ref{lm:lemma1} and assume that
(\ref{e:lmaassumption}) holds.
Let $h=h_n$ be a bandwidth and let $\beta_n = h^2 + h/\gamma_{n1}$.
Assume that $h\to0$ and $\beta_n^{-1}(\log n/n)^{1-2/\lambda} = o(1)$
For any nonnegative integer $p$, let
\[
D_{p,n}(t)={1\over n} \sum_{i=1}^n  \Biggl[{1\over m_i}
\sum_{j=1}^{m_i} K_{h,(p)} (T_{ij}-t)\mathcal{Z}_{ij} \Biggr].
\]
Then we have
\[
\sup_{t\in[a,b]}\sqrt{nh^2/(\beta_n\log n)}|D_{p,n}(t)-\mathbb
{E}\{
D_{p,n}(t)\}|=O(1)  \qquad\mbox{a.s.}
\]
\end{lm}
\begin{pf} Since both
$K$ and $t^p$ are bounded variations, $K_{(p)}$ is also a bounded variation.
Thus, we can write $K_{(p)} = K_{(p),1} - K_{(p),2}$ where $K_{(p),1}$
and $K_{(p),2}$ are
both increasing functions; without loss of generality,
assume that $K_{(p),1}(-1)=K_{(p),2}(-1)=0$.
Below, we apply Lemma~\ref{lm:lemma1} by letting $c_n=2h$.
It is clear that the assumptions of Lemma~\ref{lm:lemma1} hold here.
Write
\begin{eqnarray*}
D_n(t) &=& {1\over n}\sum_{i=1}^n  \Biggl\{{1\over m_i} \sum
_{j=1}^{m_i} K_{h,(p)}
(T_{ij}-t) \mathcal{Z}_{ij}  \Biggr\} \\
&=& {1\over n} \sum_{i=1}^n  \Biggl\{{1\over m_i} \sum
_{j=1}^{m_i}\mathcal{Z}
_{ij} I(-h\le T_{ij}-t\le h)
\int_{-h}^{T_{ij}-t} dK_{h,(p)}(v)  \Biggr\}\\
&=& \int_{-h}^{h} {1\over n} \sum_{i=1}^n  \Biggl\{{1\over m_i} \sum
_{j=1}^{m_i} \mathcal{Z}_{ij} I(v\le T_{ij}-t\le h)  \Biggr\}\,
dK_{h,(p)}(v) \\
&=& \int_{-h}^{h} G_n(t+v,t+h) \,dK_{h,(p)}(v),
\end{eqnarray*}
where $G_n$ is as defined in (\ref{e:Gn1}).
We have
%
%
\begin{eqnarray}\label{e:bddvar}
&&{\sup_{t\in[a,b]}}|D_{p,n}(t)-\mathbb{E}\{D_{p,n}(t)\}|
\nonumber\\
&&\qquad\le \sup
_{t\in
[a,b]} V_n(t,2h) \int_{-h}^{h}\bigl|dK_{h,(p)}\bigr|\\
&&\qquad\le \bigl\{K_{(p),1}(1)+ K_{(p),2}(1)\bigr\} h^{-1}
\sup_{t\in[a,b]} V_n(t,2h),\nonumber
\end{eqnarray}
and the conclusion of the lemma follows from Lemma~\ref{lm:lemma1}.
\end{pf}
\begin{pf*}{Proof of Theorem~\ref{thm:ucr}}
Define
\[
R_r^\ast= R_r -\mu(t)S_r- h\mu^{(1)}(t)S_{r+1}.
\]
By straightforward calculations, we have
%
%
\begin{equation}\label{eq:loc_lin_sol_equiv}
\widehat\mu(t)-\mu(t) = {R_0^\ast S_2-R_1^\ast S_1\over S_0S_2-S_1^2},
\end{equation}
where $S_0, S_1,S_2$ are defined as in (\ref{eq:local_linear_soluation}).
Write
\begin{eqnarray*}
R_r^\ast&=&{1\over n} \sum_{i}  \Biggl[{1\over m_i} \sum_{j=1}^{m_i}
K_h(T_{ij}-t) \{(T_{ij}-t) / h\}^r
\bigl\{Y_{ij}-\mu(t)-\mu^{(1)}(t)(T_{ij}-t)\bigr\}  \Biggr]\\
&=&{1\over n} \sum_{i}  \Biggl[{1\over m_i} \sum_{j=1}^{m_i}
K_h(T_{ij}-t) \{(T_{ij}-t) / h\}^r\\
&&\hspace*{60.3pt}{}\times
\bigl\{\varepsilon_{ij}+\mu(T_{ij})-\mu(t)-\mu^{(1)}(t)(T_{ij}-t)\bigr\}  \Biggr].
\end{eqnarray*}
By Taylor's expansion and Lemma~\ref{lm:urc_general_kernel}, uniformly
in $t$,
%
%
\begin{equation}\label{e:R*a}
R_r^\ast= {1\over n } \sum_{i} {1\over m_i} \sum_j K_h(T_{ij}-t)\{
(T_{ij}-t)/ h\}^r
\varepsilon_{ij} +O(h^2),
\end{equation}
and it follows from Lemma~\ref{lm:urc_general_kernel} that
%
%
\begin{equation}\label{e:R*b}
R_i^\ast= O\bigl(h^2 +\delta_{n1}(h)\bigr)  \qquad\mbox{a.s.}
\end{equation}
Now, at any interior point $t\in[a+h,b-h]$, since $f$ has a bounded derivative,
\begin{eqnarray*}
\mathbb{E}\{S_0\} &=& \int_{-1}^1 K(v)f(t+hv)\,dv = f(t) + O(h),\\
\mathbb{E}\{S_1\} &=& O(h),\qquad
\mathbb{E}\{S_2\}=f(t) \nu_2+O(h),
\end{eqnarray*}
where $\nu_2= \int v^2K(v)\,dv$.
By Lemma~\ref{lm:urc_general_kernel}, we conclude that, uniformly for
$t\in[a+h,b-h]$,
%
%
\begin{eqnarray}\label{e:Si}
S_0 &=& f(t) + O\bigl(h+\delta_{n1}(h)\bigr),\qquad
S_1=O\bigl(h+\delta_{n1}(h)\bigr), \nonumber\\[-8pt]\\[-8pt]
S_2 &=& f(t) \nu_2+O\bigl(h+\delta_{n1}(h)\bigr).\nonumber
\end{eqnarray}
Thus, the rate in the theorem is established by applying (\ref
{eq:loc_lin_sol_equiv}).
The same rate can also be similarly seen to hold for boundary points.
\end{pf*}

\subsection{\texorpdfstring{Proofs of Theorems \protect\ref{thm:ucr_cov} and \protect\ref{thm:error_var}}
{Proofs of Theorems 3.3 and 3.4}}
%
\begin{lm}\label{lm:ucr_cov_unif_ker}
Assume that
%
%
\begin{equation}\label{e:lmaassumption2}
\mathbb{E} \Bigl({\sup_{t\in[a,b]}} |X(t)|^{2\lambda} \Bigr) <\infty
\quad\mbox{and}\quad \mathbb{E}|U|^{2\lambda} < \infty
\qquad \mbox{for some $\lambda\in(2,\infty)$}.\hspace*{-35pt}
\end{equation}
Let $\mathcal{Z}_{ijk}$ be $X(T_{ij})X(T_{ik})$, $X(T_{ij})U_{ik}$ or
$U_{ij}U_{ik}$.
Let $c_n$ be any positive sequence tending to $0$ and
$\beta_{n} = c_n^4 + c_n^3/\gamma_{n1}+c_n^2/\gamma_{n2}$.
Assume that\break $\beta_{n}^{-1}(\log n/n)^{1-2/\lambda} = o(1)$.
Let
%
%
\begin{eqnarray}\label{e:Gn2}
&& G_n(s_1,t_1,s_2,t_2) \nonumber\\
&&\qquad= {1\over n}\sum_{i=1}^n  \biggl\{{1\over N_i}
\sum_{k\not=j}\mathcal{Z}_{ijk} I(T_{ij}\in[s_1\wedge s_2,s_1\vee
s_2],\\
&&\hspace*{133pt}T_{ik}\in[t_1\wedge t_2,t_1\vee t_2]) \biggr\}, \nonumber
\end{eqnarray}
$G(s_1,t_1,s_2,t_2)=\mathbb{E}\{G_n(s_1,t_1,s_2,t_2)\}$ and
\[
V_n(s,t,\delta)={\sup_{|u_1|,|u_2|\le\delta}}
|G_n(s,t,s+u_1,t+u_2)-G(s,t,s+u_1,t+u_2)|.
\]
Then
\[
\sup_{s,t\in[a,b]} V_n(s,t, c_n)=O(n^{-1/2} \{\beta_n \log n\}
^{1/2})  \qquad\mbox{a.s.}
\]
\end{lm}
\begin{pf} The proof is similar to that of Lemma~\ref{lm:lemma1}, and
so we only outline the main differences. 
Let $a_n,Q_n$ be as in (\ref{e:qncnan2}).
Let $\mathcal{G}$ be a two-dimensional grid on
$[a,b]^2$ with mesh $c_n$, that is, $\mathcal{G}=\{(v_{k_1},v_{k_2})\}$
where $v_k$ is defined as in the proof of
Lemma~\ref{lm:lemma1}. Then we have
%
%
\begin{equation}\label{eq:2d_grid1}
\sup_{s,t\in[a,b]} V_n(s,t, c_n) \le4
\sup_{(s,t)\in\mathcal{G}} V_n(s,t,c_n).\vadjust{\goodbreak}
\end{equation}
Define $G_n^*(s_1,t_1,s_2,t_2),G^\ast(s_1,t_1,s_2,t_2)$ and $V_n^\ast
(s,t,\delta)$
in the same way as $G_n(s_1,t_1,s_2,t_2)$, $G(s_1,t_1$, $s_2,t_2)$ and
$V_n(s,t,\delta)$
except with $\mathcal{Z}_{ijk}I(\mathcal{Z}_{ijk}\le Q_n)$ replacing
$\mathcal{Z}_{ijk}$.
Then
%
%
\begin{equation}\label{eq:2d_grid2}
\sup_{(s,t)\in\mathcal{G}} V_n(s,t,c_n)\le\sup_{(s,t)\in
\mathcal{G}}
V^\ast_n(s,t,c_n) +A_{n1}+A_{n2},
\end{equation}
where
\begin{eqnarray*}
A_{n1} &=& {\sup_{(s,t)\in\mathcal{G}} \sup_{|u_1|,|u_2|\le c_n}}
|G_n(s,t,s+u_1,t+u_2)-G_n^\ast(s,t,s+u_1,t+u_2)|, \\
A_{n2} &=& {\sup_{(s,t)\in\mathcal{G}} \sup_{|u_1|,|u_2|\le c_n}}
|G(s,t,s+u_1,t+u_2)-G^\ast(s,t,s+u_1,t+u_2)|.
\end{eqnarray*}
Using the technique similar to that in the proof of Lemma
\ref{lm:lemma1}, we can show $A_{n1}$ and $A_{n2}$ is $o(a_n)$
almost surely. To bound $V_n^\ast(s,t,c_n)$ for fixed $(s,t)$, we
create a further partition.
Put $w_n=[Q_nc_n/a_n+1]$ and $u_{r}=r c_n/w_n, r=-w_n,\ldots
, w_n$. Then
\begin{eqnarray*}
V_n^\ast(s,t,c_n)
&\le& \max_{-w_n\le r_1,r_2\le w_n} \xi_{n,r_1,
r_2}\\
&&{}+\max_{-w_n\le r_1,r_2\le w_n} \{G^\ast(s,t,s+u_{r_1+1},t+u_{r_2+1})\\
&&\hspace*{81.23pt}{}-
G^\ast(s,t,s+u_{r_1},t+u_{r_2})\},
%
\end{eqnarray*}
where
\[
\xi_{n,r_1,r_2}=|G_n^\ast(s,t,s+u_{r_1},t+u_{r_2})-G^\ast
(s,t,s+u_{r_1},t+u_{r_2})|.
\]
It is easy to see that $\operatorname{var}(\xi_{n,r_1,r_2})\le
Mn\beta_n$ for some
finite $M$,
and the rest of the proof completely mirrors that of Lemma \ref
{lm:lemma1} and is omitted.
\end{pf}
\begin{lm}\label{lm:ucr_cov_general_ker}
Let $\mathcal{Z}_{ijk}$ be as in Lemma~\ref{lm:ucr_cov_unif_ker} and
assume that
(\ref{e:lmaassumption2}) holds.
Let $h=h_n$ be a bandwidth and let $\beta_{n} = h^4 + h^3/\gamma
_{n1}+h^2/\gamma_{n2}$.
Assume that $h\to0$ and $\beta_n^{-1}(\log n/n)^{1-2/\lambda} = o(1)$.
For any nonnegative integers $p,q$, let
\[
D_{p,q,n}(s,t)={1\over n }\sum_{i=1}^n  \biggl[ {1\over N_j}
\sum_{k\not=j} \mathcal{Z}_{ijk} K_{h,(p)}(T_{ij}-s)
K_{h,(q)}(T_{ik}-t) \biggr].
\]
Then, for any $p,q$,
\[
\sup_{s,t\in[a,b]} \sqrt{nh^4/(\beta_n\log n)}
|D_{p,q,n}(s,t)-\mathbb{E}
\{D_{p,q,n}(s,t)\}|=O(1)  \qquad\mbox{a.s.}
\]
\end{lm}
\begin{pf} Write
{\fontsize{10.5}{11}\selectfont{\begin{eqnarray*}
\hspace*{-4pt}&&
D_{p,q,n}(s,t)\\
\hspace*{-4pt}&&\qquad=\sum_{i=1}^n
\biggl[ {1\over N_i} \sum_{k\not=j} \mathcal{Z}_{ijk} I(T_{ij}\le s+h)
I(T_{ik}\le t+h) \\
\hspace*{-4pt}&&\qquad\quad\hspace*{51pt}{} \times
K_{h,(p)}(T_{ij}-s) K_{h,(q)}(T_{ik}-t)  \biggr]\\
\hspace*{-4pt}&&\qquad= \int\int_{(u,v)\in[-h,h]^2} {1\over n }\sum_{i=1}^n
\biggl[ {1\over N_i} \sum_{k\not=j} \mathcal{Z}_{ijk}\\
\hspace*{-4pt}&&\qquad\quad\hspace*{121.7pt}{}\times I(T_{ij}\in
[s+u,s+h])\\
\hspace*{-4pt}&&\qquad\quad\hspace*{122.5pt}{} \times I(T_{ik}\in[t+v,t+h]) \biggr]
\,d K_{h,(p)}(u)\,dK_{h,(q)}(v)\\
\hspace*{-4pt}&&\qquad= \int\int_{(u,v)\in[-h,h]^2} G_n(s+u,t+v,s+h,t+h) \,d K_{h,(p)}(u)
d K_{h,(q)}(v),
\end{eqnarray*}}}

\noindent where $G_n$ is as in (\ref{e:Gn2}). Now,
\begin{eqnarray*}
&&{\sup_{(s,t)\in[a,b]^2}}|D_{p,q,n}(s,t)-\mathbb{E}\{D_{p,q,n}(s,t)\}
|\\
&&\qquad\le
\sup_{s,t\in[a,b]} V_n(s,t,2h)
\int\int_{(u,v)\in[-h,h]^2} \bigl|d\bigl\{K_{h,(p)}(u)\bigr\}\bigr| \bigl|d\bigl\{
K_{h,(q)}(v)\bigr\}\bigr|\\
&&\qquad= O[\{\beta_n\log n/(nh^4)\}^{1/2}] \qquad\mbox{a.s.}
\end{eqnarray*}
by Lemma~\ref{lm:ucr_cov_unif_ker}, using the same argument as in
(\ref{e:bddvar}).
\end{pf}
\begin{pf*}{Proof of Theorem~\ref{thm:ucr_cov}}
Let $S_{pq}, R_{pq}, \mathcal{A}_i$ and $\mathcal{B}$ be defined as
in (\ref
{eq:cov_est_sol}).
Also, for $p,q\ge0$, define
\[
R_{pq}^\ast= R_{pq} - C(s,t) S_{pq}- h_RC^{(1,0)}(s,t) S_{p+1,q} -
h_RC^{(0,1)}(s,t) S_{p,q+1}.
\]
By straightforward algebra, we have
%
%
\begin{equation}\label{e:hatC-C}
(\widehat C-C)(s,t)=(\mathcal{A}_1 R_{00}^\ast-\mathcal{A}_2
R_{10}^\ast-\mathcal{A}_3
R_{01}^\ast)\mathcal{B}^{-1}.
\end{equation}
By standard calculations, we have
the following rates uniformly on $[a+h_R,b-h_R]^2$:
\begin{eqnarray*}
\mathbb{E}(S_{00})&=&f(s)f(t)+O(h_R),\qquad
\mathbb{E}(S_{01})=O(h_R),\\
\mathbb{E}(S_{10})&=&O(h_R),\qquad
\mathbb{E}(S_{02})=f(s)f(t)\nu_2 +O(h_R),\\
\mathbb{E}(S_{20})&=&f(s)f(t)\nu_2 +O(h_R),\qquad
\mathbb{E}(S_{11})=O(h_R).
%
%
\end{eqnarray*}
By these and Lemma
\ref{lm:ucr_cov_general_ker}, we have the following almost sure
uniform rates:
%
%
\begin{eqnarray}\label{e:AB}
\mathcal{A}_1 &=& f^2(s) f^2(t)\nu_2^2 +O\bigl(h_R+\delta_{n2}(h_R)\bigr), \nonumber\\
\mathcal{A}_2 &=&
O\bigl(h_R+\delta_{n2}(h_R)\bigr),\nonumber\\[-8pt]\\[-8pt]
\mathcal{A}_3 &=& O\bigl(h_R+\delta_{n2}(h_R)\bigr), \nonumber\\
\mathcal{B}&=& f^3(s)f^3(t)\nu_2^2
+O\bigl(h_R+\delta_{n2}(h_R)\bigr).\nonumber
\end{eqnarray}
To analyze the behavior of the components of (\ref{e:hatC-C}), it
suffices now to
analyze $R_{pq}^\ast$. Write
\begin{eqnarray*}
R_{00}^\ast
&=& {1\over n} \sum_{i=1}^n  \biggl[ {1\over N_i}\sum_{k\not=j}
\bigl\{Y_{ij} Y_{ik}-C(s,t)\\
&&\hspace*{64.35pt}{}-C^{(1,0)}(s,t) (T_{ij}-s)  \\
&&\hspace*{64.35pt}{}- C^{(0,1)}(s,t)(T_{ik}-t)\bigr\}\times K_{h_R}(T_{ij}-s)
K_{h_R}(T_{ik}-t) \biggr].
\end{eqnarray*}
Let $\varepsilon_{ijk}^\ast=Y_{ij}Y_{ik}-C(T_{ij},T_{ik})$.
By Taylor's expansion,
\begin{eqnarray*}
&&Y_{ij} Y_{ik}-C(s,t)-C^{(1,0)}(s,t) (T_{ij}-s)
-C^{(0,1)}(s,t)(T_{ik}-t) \\
&&\qquad= Y_{ij} Y_{ik} - C(s,t) - C(T_{ij},T_{ik}) + C(T_{ij},T_{ik})
\\
&&\qquad\quad\hspace*{0pt}{}-C^{(1,0)}(s,t) (T_{ij}-s) -C^{(0,1)}(s,t)(T_{ik}-t) \\
&&\qquad= \varepsilon_{ijk}^\ast+ O(h_R^2) \qquad \mbox{a.s.}
\end{eqnarray*}
It follows that
%
%
\begin{equation}\label{e:R00*}
R_{00}^\ast= {1\over n}\sum_{i=1}^n {1\over N_i}\sum_{k\not=j}
\varepsilon_{ijk}^\ast
K_{h_R}(T_{ij}-s) K_{h_R}(T_{ik}-t) + O(h_R^2).
\end{equation}
Applying Lemma~\ref{lm:ucr_cov_general_ker}, we obtain, uniformly in $s,t$,
%
%
\begin{equation}\label{e:R00*a}
R_{00}^\ast= O\bigl(\delta_{n2}(h_R)+h_R^2\bigr)\qquad  \mbox{a.s.}
\end{equation}
By (\ref{e:AB}),
%
%
\begin{equation}\label{e:ABinv}
\mathcal{A}_1 \mathcal{B}^{-1} = [f(s)f(t)]^{-1} + O\bigl(h_R+\delta_{n2}(h_R)\bigr).
\end{equation}
Thus, $R_{00}^\ast\mathcal{A}_1\mathcal{B}^{-1} = O(\delta
_{n2}(h_R)+h_R^2)$ a.s.
Similar derivations show that\break $R_{10}^\ast\mathcal{A}_2\times\mathcal
{B}^{-1}$ and
$R_{01}^\ast\mathcal{A}_3\mathcal{B}^{-1}$ are both of lower order.
Thus, the rate in (\ref{e:Rrate}) is obtained for $s,t\in[a+h_R,b-h_R]$.
As for $s$ and/or $t$ in $[a,a+h)\cup(b-h,b]$,
similar calculations show that the same rate also holds. The result follows
by taking into account of the rate of $\widehat\mu$.
\end{pf*}
\begin{pf*}{Proof of Theorem~\ref{thm:error_var}}
Note that
\[
\widehat\sigma^2 -\sigma^2 = {1\over b-a} \int_a^b \{\widehat
V(t)- V(t)\} \,dt
- {1\over b-a} \int_a^b \{\widehat C(t,t)-C(t,t)\} \,dt.
\]
To consider $\widehat V(t)- V(t)$ we follow the development in the
proof of
Theorem~\ref{thm:ucr}.
Recall (\ref{eq:lls_error-var}) and let $Q_r^\ast= Q_r -V(t)S_r- h
V^{(1)}(t)S_{r+1}$.
Then, as in~(\ref{eq:loc_lin_sol_equiv}), we obtain
\[
\widehat V(t)- V(t) = {Q_0^\ast S_2-Q_1^\ast S_1\over S_0S_2-S_1^2}.
\]
Write
\begin{eqnarray*}
Q_r^\ast&=&{1\over n} \sum_{i=1}^n
{1\over m_i} \sum_{j=1}^{m_i} K_{h_V}(T_{ij}-t) \{(T_{ij}-t) / h_V\}^r
\bigl\{Y_{ij}^2-V(t)-V^{(1)}(t)(T_{ij}-t)\bigr\} \\[-2pt]
&=&{1\over n} \sum_{i=1}^n {1\over m_i} \sum_{j=1}^{m_i}
K_{h_V}(T_{ij}-t) \{(T_{ij}-t) / h\}^r
\{Y_{ij}^2-V(T_{ij})\} + O(h_V^2),
\end{eqnarray*}
which, by Lemma~\ref{lm:lemma1}, has the uniformly rate $O(h_V^2 +
\delta_{n1}(h_V))$ a.s.
By (\ref{e:Si}), we have
\begin{eqnarray*}
\widehat V(t)- V(t)
&=& {1\over f(t)n} \sum_{i} {1\over m_i} \sum_{j=1}^{m_i}
K_h(T_{ij}-t)
\{Y_{ij}^2-V(T_{ij})\}\\[-2pt]
&&{} + O\bigl(h_V^2+\delta_{n1}^2(h_V)\bigr)  \qquad\mbox{a.s.}
\end{eqnarray*}
Thus,
\begin{eqnarray*}
\int_a^b \{\widehat V(t)- V(t)\} \,dt
&=& {1\over n} \sum_{i=1}^n {1\over m_i} \sum_{j=1}^{m_i}
\{Y_{ij}^2-V(T_{ij})\} \int_a^b K_{h_V}(T_{ij}-t) f^{-1}(t)\,dt \\[-2pt]
&&{} + O\bigl(h_V^2+\delta_{n1}^2(h_V)\bigr)\qquad
\mbox{a.s.}
\end{eqnarray*}
Note that
\[
\biggl|\int_a^b K_{h_V}(T_{ij}-t) f^{-1}(t)\,dt \biggr| \le\sup_t f^{-1}(t).
\]
By Lemma~\ref{lm:root-n} below in this subsection,
%
%
\begin{equation}\label{e:sigma1}\quad
\int_a^b \{\widehat V(t)- V(t)\} \,dt = O\bigl((\log n/n)^{1/2}+h_V^2+\delta
_{n1}^2(h_V)\bigr)\qquad
\mbox{a.s.}
\end{equation}

Next, we consider $\widehat C(t,t)-C(t,t)$. We apply (\ref{e:hatC-C}) but
will focus on
$R_{00}^\ast\mathcal{A}_1 \mathcal{B}^{-1}$ since the other two
terms are dealt with
similarly.
By (\ref{e:R00*})--(\ref{e:ABinv}),
%
%
\begin{eqnarray}\label{e:R00*a}\quad
R_{00}^\ast\mathcal{A}_1 \mathcal{B}^{-1} &=& {1\over f(s)f(t)n}\sum_{i=1}^n
{1\over N_i}
\sum_{k\not=j} \varepsilon_{ijk}^\ast K_{h_R}(T_{ij}-s)
K_{h_R}(T_{ik}-t) \nonumber\\[-9pt]\\[-9pt]
&&{} + O\bigl(h_R^2+\delta_{n2}^2(h_R)\bigr)  \qquad\mbox{a.s.}
\nonumber
\end{eqnarray}
Thus,
\begin{eqnarray*}
&&\int_a^b \{\widehat C(t,t)-C(t,t)\} \,dt
\\[-2pt]
&&\qquad= {1\over n}\sum_{i=1}^n {1\over N_i}
\sum_{k\not=j} \varepsilon_{ijk}^\ast\int_a^b K_{h_R}(T_{ij}-t)
K_{h_R}(T_{ik}-t) f^{-2}(t) \,dt \\[-2pt]
&&\qquad\quad{} + O\bigl(h_R^2+\delta_{n2}^2(h_R)\bigr)  \qquad\mbox{a.s.}
\end{eqnarray*}
Write
\begin{eqnarray*}
&&
\int_a^b K_{h_R}(T_{ij}-t) K_{h_R}(T_{ik}-t) f^{-2}(t) \,dt
\\
&&\qquad= \int_{-1}^1 K(u)K_{h_R}\bigl((T_{ik}-T_{ij})+u h_R\bigr)f^{-2}(T_{ij}- u h_R)\,du.
\end{eqnarray*}
%
%
%
%
%
A slightly modified version of Lemma~\ref{lm:lemma1} leads to the
``one-dimensional'' rate:
\begin{eqnarray*}
&&\sup_{u\in[0,1]}  \Biggl|{1\over n}\sum_{i=1}^n {1\over N_i}\sum
_{k\not=j} \varepsilon_{ijk}^\ast
K_{h_R}\bigl((T_{ik}-T_{ij})+uh_R\bigr)f^{-2}(T_{ij}-uh_R) \Biggr|
\\
&&\qquad= O(\delta_{n1}(h_R))  \qquad\mbox{a.s.}
\end{eqnarray*}
It follows that
%
%
\begin{equation}\label{e:sigma2}\quad
\int_a^b \{\widehat C(t,t)-C(t,t)\} \,dt = O\bigl(h_R^2+\delta
_{n1}(h_R)+\delta
_{n2}^2(h_R)\bigr)  \qquad\mbox{a.s.}
\end{equation}
The theorem follows from (\ref{e:sigma1}) and (\ref{e:sigma2}).
\end{pf*}
\begin{lm}\label{lm:root-n}
Assume that $\xi_{ni},1\le i\le n$, are independent random variables
with mean zero and
finite variance. Also assume that there exist i.i.d. random variables
$\xi_i$ with mean
zero and finite $\delta$th moment for some $\delta>2$ such that $|\xi
_{ni}|\le|\xi_i|$. Then
\[
{1\over n} \sum_{i=1}^n\xi_{ni} = O\bigl((\log n/n)^{1/2}\bigr) \qquad \mbox{a.s.}
\]
\end{lm}
\begin{pf}
Let $a_n = (\log n/n)^{1/2}$. Assume that $\xi_{ni}\ge0$.
Write
\[
\xi_{ni} = \xi_{ni\succ}+\xi_{ni\prec} := \xi_{ni}I(|\xi_{ni}|> a_n^{-1})
+\xi_{ni}I(|\xi_{ni}|\le a_n^{-1}).
\]
Then
\[
\Biggl|{1\over a_nn} \sum_{i=1}^n\xi_{ni\succ} \Biggr|
\le{1\over a_nn} \sum_{i=1}^n |\xi_{ni\succ}|^{\delta}
|\xi_{ni\succ}|^{1-\delta} \le a_n^{\delta-2} {1\over n} \sum
_{i=1}^n|\xi_{i}|^{\delta}\to0
\qquad\mbox{a.s.}
\]
by the law of large numbers. The mean of the left-hand side is also
tending to zero
by the same argument. Thus, $n^{-1} \sum_{i=1}^n(\xi_{ni\succ
}-\mathbb{E}\{
\xi_{ni\succ}\})
=o(a_n)$. Next, by Bernstein's inequality,
\begin{eqnarray*}
\mathbb{P} \Biggl({1\over n} \sum_{i=1}^n(\xi_{ni\prec}-\mathbb{E}\{
\xi_{ni\prec
}\}) > Ba_n \Biggr)
&\le&\exp \biggl\{-{B^2n^2a_n^2\over2 n \sigma^2 + (2/3)Bn} \biggr\}
\\
&=& \exp \biggl\{-{B^2\log n\over2\sigma^2 + (2/3)B} \biggr\},
\end{eqnarray*}
which is summable for large enough $B$. The result follows from the
Borel--Cantelli lemma.
\end{pf}

%
\subsection{\texorpdfstring{Proof of Theorem \protect\ref{thm:pca_rate}}{Proof of Theorem 3.6}}
\label{sec:proof_thm_dev_psi}

Let $\Delta$ be the integral operator with kernel $\widehat R-R$.
\begin{lm} \label{lm:Delta}
For any bounded measurable function $\psi$ on $[a,b]$,
\[
{\sup_{t\in[a,b]}} |(\Delta\psi)(t)| = O\bigl(h_{\mu}^2+\delta
_{n1}(h_\mu)+h_R^2
+\delta_{n1}(h_R)+\delta_{n2}^2(h_R)\bigr) \qquad\mbox{a.s.}
\]
\end{lm}
\begin{pf} It follows that
\begin{eqnarray*}
(\Delta\psi)(t) &=& \int_{s=a}^b (\widehat C-C)(s,t) \psi(s) \,ds
- \int_{s=a}^b \{\widehat\mu(s)\widehat\mu(t)-\mu(s)\mu(t)\}
\psi(s)\,ds \\
&=&\!:
A_{n1}-A_{n2}.
\end{eqnarray*}
By (\ref{e:hatC-C}),
\[
A_{n1}=\int_{s=a}^b (\mathcal{A}_1 R_{00}^\ast-\mathcal{A}_2
R_{10}^\ast-\mathcal{A}_3
R_{01}^\ast)\mathcal{B}^{-1} \psi(s) \,ds.
\]
We focus on $\int_{s=a}^b \mathcal{A}_1 R_{00}^\ast\mathcal{B}^{-1}
\psi(s) \,ds$
since the other two terms are of
lower order and can be dealt with similarly. By (\ref{e:R00*}) and
(\ref{e:ABinv}),
\begin{eqnarray*}
&& \int_{s=a}^b \mathcal{A}_1 R_{00}^\ast\mathcal{B}^{-1} \psi
(s)\,ds \\
&&\qquad= {1\over f(t)n}\sum_{i=1}^n {1\over N_i} \sum_{k\not=j} \varepsilon
_{ijk}^\ast
K_{h_R}(T_{ik}-t) \int_{s=a}^b K_{h_R}(T_{ij}-s) \psi(s) f(s)^{-1} \,ds
\\
&&\qquad\quad{} + O\bigl(h_R^2+\delta_{n2}^2(h_R)\bigr).
\end{eqnarray*}
%
%
%
%
%
%
Note that
\[
\biggl|\int_{s=a}^b K_{h_R}(T_{ij}-s) \psi(s) f(s)^{-1} \,ds \biggr|
\le\sup_{s\in[a,b]} (|\psi(s)| f(s)^{-1}) \int_{u=-1}^1 K(u)\, du.
\]
Thus, Lemma~\ref{lm:lemma1} can be easily improvised to give the
following uniform rate
over~$t$:
\begin{eqnarray*}
&&{1\over f(t)n}\sum_{i=1}^n {1\over N_i} \sum_{k\not=j} \varepsilon
_{ijk}^\ast
K_{h_R}(T_{ik}-t) \int_{s=a}^b K_{h_R}(T_{ij}-s) \psi(s) f(s)^{-1} \,ds
\\
&&\qquad= O(\delta_{n1}(h_R))\qquad\mbox{a.s.}
\end{eqnarray*}
Thus,
\[
\int_{s=a}^b \mathcal{A}_1 R_{00}^\ast\mathcal{B}^{-1} \psi(s)\,ds
= O\bigl(\delta_{n1}(h_R)+h_R^2+\delta_{n2}^2(h_R)\bigr) \qquad\mbox{a.s.},
\]
which is also the rate of $A_{n1}$.
Next, we write
\[
A_{n2}=\widehat\mu(t)\int_{s=a}^b \{\widehat\mu(s)-\mu(s)\}\psi(s)\,ds
-\{\widehat\mu(t)-\mu(t)\}\int_{s=a}^b \mu(s)\psi(s)\,ds,
\]
which has the rate $O(h_{\mu}^2+\delta_{n1}(h_\mu))$ by
Theorem~\ref{thm:ucr}.\vspace*{-3pt}
\end{pf}
\begin{pf*}{Proof of Theorem~\ref{thm:pca_rate}}
We prove (b) first. \citet{HH06} give the $L^2$ expansion
\[
\widehat\psi_j-\psi_j = \sum_{k\not=j}(\lambda_j-\lambda
_k)^{-1} {\langle}
\Delta\psi_j,\psi_k{\rangle}\phi_k
+ O(\|\Delta\|^2),
\]
where $\|\Delta\| = (\iint\{\widehat R(s,t)-R(s,t)\}^2\,ds\,dt)^{1/2}$, the
Hilbert--Schmidt norm of~$\Delta$.
By Bessel's inequality, this leads to
\[
\|\widehat\psi_j-\psi_j\| \le C(\|\Delta\psi_j\| + \|\Delta\|^2).
\]
By Lemma~\ref{lm:Delta} and Theorem~\ref{thm:ucr_cov},
\begin{eqnarray*}
\|\Delta\psi_j\| &=& O\bigl(h_{\mu}^2+\delta_{n1}(h_\mu)+h_R^2
+\delta_{n1}(h_R)+\delta_{n2}^2(h_R)\bigr), \\[-3pt]
\|\Delta\|^2 &=& O\bigl(h_\mu^4+\delta_{n1}^2(h_\mu)+ h_R^4 + \delta
_{n2}^2(h_R)\bigr)\qquad \mbox{a.s.}
\end{eqnarray*}
%
%
Thus,
\[
\|\widehat\psi_j-\psi_j\| = O\bigl(h_{\mu}^2+\delta_{n1}(h_\mu)+h_R^2
+\delta_{n1}(h_R)+\delta_{n2}^2(h_R)\bigr)\qquad \mbox{a.s.},
\]
proving (b).

Next, we consider (a). By (4.9) in \citet{HMW06},
\begin{eqnarray*}
\widehat\omega_j-\omega_j
&=& \iint(\widehat R-R)(s,t) \psi_j(s) \psi_j(t) \,ds \,dt
+ O(\|\Delta\psi_j\|^2)\\[-3pt]
&=& \iint(\widehat C-C)(s,t) \psi_j(s)\psi_j(t)\, ds\,dt \\[-3pt]
&& - \iint\{\widehat\mu(s)\widehat\mu(t)-\mu(s)\mu
(t)\}\psi
_j(s)\psi_j(t)\,ds\,dt
+ O(\|\Delta\psi_j\|^2) \\[-3pt]
&=&\!: A_{n1}-A_{n2}+ O(\|\Delta\psi_j\|^2).
\end{eqnarray*}
Now,
\[
A_{n1}=\iint(\mathcal{A}_1 R_{00}^\ast-\mathcal{A}_2 R_{10}^\ast
-\mathcal{A}_3
R_{01}^\ast)\mathcal{B}^{-1}
\psi_j(s)\psi_j(t) \,ds\,dt.
\]
Again it suffices to focus on $\iint\mathcal{A}_1 R_{00}^\ast
\mathcal{B}^{-1} \psi
_j(s)\psi_j(t) \,ds\,dt$.
By (\ref{e:R00*}) and~(\ref{e:ABinv}),
\begin{eqnarray*}
&& \iint\mathcal{A}_1 R_{00}^\ast\mathcal{B}^{-1} \psi_j(s)\psi
_j(t)\,ds\,dt \\[-3pt]
&&\qquad= {1\over n}\sum_{i=1}^n {1\over m_i(m_i-1)} \sum_{k\not=j}
\varepsilon_{ijk}^\ast
\iint K_{h_R}(T_{ij}-s)K_{h_R}(T_{ik}-t)\\[-3pt]
&&\qquad\quad\hspace*{131.4pt}{}\times \psi_j(s) \psi_j(t) \{
f(s)f(t)\}^{-1} \,ds\,dt \\[-3pt]
&&\qquad\quad{} + O\bigl(h_R^2+\delta_{n2}^2(h_R)\bigr)
\qquad\mbox{a.s.},
\end{eqnarray*}
where the first term on the right-hand side can be shown to be $O((\log
/n)^{1/2})$ a.s.
by Lemma~\ref{lm:root-n}.
Thus,
\[
A_{n1}=O\bigl((\log/n)^{1/2}+h_R^2+\delta_{n2}^2(h_R)\bigr).
\]
Next, write
\begin{eqnarray*}
A_{n2}&=& \int\{\widehat\mu(s)-\mu(s)\}\psi_j(s)\,ds
\int\widehat\mu(t)\psi_j(t)\,dt \\
&&{} + \int\mu(s)\psi_j(s)\,ds \int\{\widehat\mu(t)-\mu(t)\}\psi_j(t)\,dt,
\end{eqnarray*}
and it can be similarly shown that
\[
A_{n2} = O\bigl((\log/n)^{1/2}+h_\mu^2+\delta_{n1}^2(h_\mu)\bigr)\qquad\mbox{a.s.}
\]
This establishes (a).

Finally, we consider (c). For any $t\in[a,b]$,
\begin{eqnarray*}
&&
\widehat\omega_j \widehat\psi_j(t)-\omega_j \psi_j(t)\\
&&\qquad=\int
\widehat R(s,t) \widehat
\psi_j(s) \,ds -\int R(s,t) \psi_j(s) \,ds\\
%
%
&&\qquad= \int\{\widehat R(s,t)-R(s,t)\} \psi_j(s) \,ds + \int\widehat
R(s,t) \{\widehat
\psi_j(s)-\psi_j(s)\} \,ds.
%
%
\end{eqnarray*}
By the Cauchy--Schwarz inequality, uniformly for all $t\in[a,b]$,
\begin{eqnarray*}
\biggl|\int\widehat R(s,t) \{\widehat\psi_j(s)-\psi_j(s)\} \,ds
\biggr| &\le&
\biggl\{\int\widehat R^2(s,t) \,ds  \biggr\}^{1/2} \|\widehat\psi
_j-\psi_j\|\\
&\le& |b-a|^{1/2} \sup_{s,t} |\widehat R(s,t)| \times\|\widehat
\psi_j-\psi
_j\|\\
&=& O(\|\widehat\psi_j-\psi_j\|) \qquad\mbox{a.s.}
\end{eqnarray*}
Thus,
\[
\widehat\omega_j \widehat\psi_j(t)-\omega_j \psi_j(t) =
O\bigl(h_{\mu}^2+\delta
_{n1}(h_\mu)+h_R^2
+\delta_{n1}(h_R)+\delta_{n2}^2(h_R)\bigr)\qquad \mbox{a.s.}
\]
By the triangle inequality and (b),
\begin{eqnarray*}
&&\omega_j |\widehat\psi_j(t)-\psi_j(t)| \\
&&\qquad= | \widehat\omega
_j\widehat\psi
_j(t)-\omega_j \psi_j(t) - (\widehat\omega_j-\omega_j)\widehat
\psi_j(t) |\\
&&\qquad\le | \widehat\omega_j\widehat\psi_j(t)-\omega_j \psi_j(t)|
+ {|\widehat
\omega_j-\omega_j| \sup_t}|\widehat\psi_j(t) |\\
&&\qquad= O\bigl((\log n/n)^{1/2}+h_{\mu}^2+\delta_{n1}(h_\mu)+h_R^2
+\delta_{n1}(h_R)+\delta_{n2}^2(h_R)\bigr) \qquad\mbox{a.s.}
\end{eqnarray*}
Note that $(\log n/n)^{1/2}=o(\delta_{n1}(h_\mu))$.
This completes the proof of (c).
\end{pf*}

\section*{Acknowledgments}
We are very grateful to the Associate Editor and two referees for their
helpful comments
and suggestions.

\printaddresses

\end{document}